\documentclass[11pt,reqno]{amsart}
\usepackage{amssymb, amsmath, amsthm, amsfonts,xcolor, enumerate}
\usepackage{verbatim}

\usepackage{amsmath, amsopn, amssymb, amsthm}
\usepackage{enumerate}
\usepackage[font={small}]{caption}

\usepackage[tmargin=1.5in,bmargin=1.5in,lmargin=1.3in,rmargin=1.3in]{geometry}
\usepackage{bbm}
\usepackage{cite}
\usepackage{microtype}

\definecolor{halfgray}{gray}{0.55}
\definecolor{webgreen}{RGB}{0, 100, 0}
\definecolor{webbrown}{rgb}{.6,0,0}
\definecolor{Maroon}{cmyk}{0, 0.87, 0.68, 0.32}
\definecolor{RoyalBlue}{RGB}{0, 0, 128}
\definecolor{Black}{cmyk}{0, 0, 0, 0}

\usepackage{hyperref}
\hypersetup{%
    colorlinks=true, linktocpage=true, pdfstartpage=1, pdfstartview=FitV,%
    breaklinks=true, pdfpagemode=UseNone, pageanchor=true, pdfpagemode=UseOutlines,%
    plainpages=false, bookmarksnumbered, bookmarksopen=true, bookmarksopenlevel=1,%
    hypertexnames=true, pdfhighlight=/O,
    urlcolor=webbrown, linkcolor=RoyalBlue, citecolor=webgreen, 
}

\usepackage[textsize=small]{todonotes}       

\newtheorem{thm}{Theorem}[section]

\newtheorem{lemma}[thm]{Lemma}
\newtheorem{prop}[thm]{Proposition}

\newtheorem{defn}[thm]{Definition}

\theoremstyle{definition}

\def\A{\mathcal{A}}
\def\B{\mathcal{B}}

\def\E{\mathcal{E}}
\def\F{\mathcal{F}}

\def\HH{\mathcal{H}}

\def\M{\mathcal{M}}

\def\P{\mathcal{P}}
\def\Q{\mathcal{Q}}
\def\R{\mathcal{R}}
\def\S{\mathcal{S}}

\def\U{\mathcal{U}}

\newcommand{\Acal}{\mathcal{A}}
\newcommand{\Bcal}{\mathcal{B}}

\newcommand{\Ncal}{\mathcal{N}}

\newcommand{\Pp}{\mathbb{P}_p}
\newcommand{\Qcal}{\mathcal{Q}}

\newcommand{\V}{\mathcal{V}}

\def\Ex{\mathbb{E}}

\def\N{\mathbb{N}}
\def\Pr{\mathbb{P}}

\def\RR{\mathbb{R}}
\def\ZZ{\mathbb{Z}}
\newcommand{\Hb}{\mathbb{H}}

\def\le{\leqslant}
\def\ge{\geqslant}

\def\->{\rightarrow}
\def\<{\langle}
\def\>{\rangle}

\newcommand{\growto}{\Rightarrow}

\newcommand{\vep}{\varepsilon}

\newcommand{\sss}{\scriptscriptstyle}
\newcommand{\e}{\mathrm{e}}
\renewcommand{\d}{\mathrm{d}}
\newcommand{\sq}{\succcurlyeq}

\newcommand{\ul}{\underline}
\renewcommand{\emptyset}{\varnothing}
\newcommand{\indi}{\mathbbm{1}}
\newcommand{\Mbb}{\mathbb{M}}
\newcommand{\mfr}{\mathfrak{r}}
\newcommand{\mfR}{\mathfrak{R}}
\newcommand{\xb}{\mathbf{x}}
\newcommand{\yb}{\mathbf{y}}

\def\le{\leqslant}
\def\ge{\geqslant}

\def\Nani{\mathcal{N}_{\sss (1,2)}}
\def\Ngen{\mathcal{N}_{\sss (a,b)}}
\def\Nnn{\mathcal{N}_{\sss (1,1)}}
\def\Nnnd{\mathcal{N}_{\sss (1,\ldots,1)}}
\def\Ndu{\mathcal{N}_{\textup{Duarte}}}

\def\<{\langle}
\def\>{\rangle}

\numberwithin{equation}{section}

\title{Higher order corrections for anisotropic bootstrap percolation}

\begin{document}
\author[H. Duminil-Copin]{Hugo Duminil-Copin}
\address{Institut des Hautes \'Etudes Scientifiques, Le Bois-Marie 35, route de Chartres 91440 Bures-sur-Yvette France, and D\'epartement de math\'ematiques Ð Universit\'e de Gen\`eve, 2-4 rue du li\`evre, 64 1211 Gen\`eve 4, Switzerland}
\email{duminil@ihes.fr}

\author[A.C.D. van Enter]{Aernout C.D. van Enter}
\address{Johann Bernoulli Instituut voor Wiskunde en Informatica, University of Groningen, PO Box 407, 9700 AK, Groningen, the Netherlands}
\email{avanenter@gmail.com}

\author[T. Hulshof]{Tim Hulshof}
\address{Department of Mathematics and Computer Science, Eindhoven University of Technology, PO Box 513, 5600 MB Eindhoven, the Netherlands}\email{w.j.t.hulshof@tue.nl}

\begin{abstract}
We study the critical probability for the metastable phase transition of the two-dimensional anisotropic bootstrap percolation model with $(1,2)$-neighbourhood and threshold $r = 3$. The first order asymptotics for the critical probability were recently determined by the first and second authors. Here we determine the following sharp second and third order asymptotics:
\begin{multline*}
	p_c\big( [L]^2,\Nani,3 \big) \; = \;  \frac{(\log \log L)^2}{12\log L} \, - \, \frac{\log \log L \, \log \log \log L}{ 3\log L}
	 +  \frac{\left(\log \frac{9}{2}  + 1 \pm o(1) \right)\log \log L}{6\log L}.
\end{multline*}
We note that the second and third order terms are so large that the first order asymptotics fail to approximate $p_c$ even for lattices of size well beyond $10^{10^{1000}}$.
\end{abstract}

\date{\today}
\maketitle
\vspace{1em}
{\small
\noindent
{\it MSC 2010.} 60K35, 82B43, 82C43.

\noindent
{\it Key words and phrases.}
Bootstrap percolation, finite-size effects, metastability, sharp threshold. 
}
\vspace{1em}
\hrule

\section{Introduction}\label{sec:intro}

\subsection{Motivation and statement of the main result}

Bootstrap percolation is a general name for the dynamics of monotone, two-state cellular automata on a graph $G$.
Bootstrap percolation models with different rules and on different graphs have since their invention by Chalupa, Leath and Reich \cite{ChaLeaRei79} been applied in various contexts 
and the mathematical properties of bootstrap percolation are an active area of research at the intersection between probability theory and combinatorics. See for instance 
\cite{AizLeb88,Adl91, AdlLev03, Ami10, Baletal12, Hol03, DumEnt13} and the references therein.

Motivated by applications to statistical (solid-state) physics such as the Glauber dynamics of the Ising model~\cite{FonSchSid02,Mor11} and kinetically constrained spin models~\cite{Canetal08}, the underlying graph is often taken to be a $d$-dimensional lattice, and the initial state is usually chosen randomly. 

Although some progress has recently been made in the study of very general cellular automata on lattices~\cite{BolSmiUzz12a,Boletal14,DumHol12}, attention so far has mainly focused on obtaining a very precise understanding of the metastable transition for specific simple models~\cite{AizLeb88,Hol03,BalBolMor09,Baletal12,CerCir99,DumEnt13,GraHolMor12}. 

In this paper we will provide the most detailed description so far for such a model; namely, the so-called \emph{anisotropic bootstrap percolation model}, defined as follows: First, given a finite set $\Ncal \subset \ZZ^d \setminus \{0\}$ (the \emph{neighbourhood}) and an integer $r$ (the \emph{threshold}), define the bootstrap operator
\begin{equation}\label{def:operator}
	\B(\S) \, := \, \S \cup \big\{ v \in \ZZ^d \,:\, | (v + \Ncal) \cap \S | \ge r \big\}
\end{equation}
for every set $\S \subset \ZZ^d$. That is, viewing $\S$ as the set of ``infected'' sites, every site $v$ that has at least $r$ infected ``neighbours'' in $v+\Ncal$ becomes infected by the application of $\B$. For $t \in \mathbb{N}$ let $\B^{(t)} (\S) = \B (\B^{(t-1)}(\S))$, where $\Bcal^{(0)}(\S) = \S$, and let $\< \S \> = \lim_{t \to \infty} \mathcal{B}^{(t)}(\S)$ denote the set of eventually infected sites. For each $p \in (0,1)$, let $\Pp$ denote the probability measure under which the elements of the initial set $\S \subset \ZZ^d$ are chosen independently at random with probability $p$, and for each set $\Lambda \subset \ZZ^d$, define the \emph{critical probability} on $\Lambda$ to be
\begin{equation}\label{def:pc}
p_c\big( \Lambda, \Ncal, r \big) \, := \, \inf\big\{ p > 0 \,:\, \Pp \big( \Lambda \subset \< \S \cap \Lambda \> \big) \ge 1/2 \big\}.
\end{equation}
If $\Lambda \subset \< \S \cap \Lambda \>$ then we say that $\S$ \emph{percolates} on $\Lambda$. We remark that since we will usually expect the probability of percolation to undergo a sharp transition around $p_c$, the choice of the constant $1/2$ in the definition~\eqref{def:pc} is not significant. 

\begin{figure}
	\includegraphics[width = .49\textwidth]{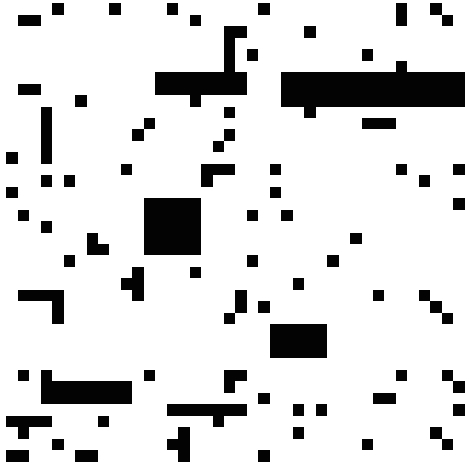}
	\includegraphics[width = .49\textwidth]{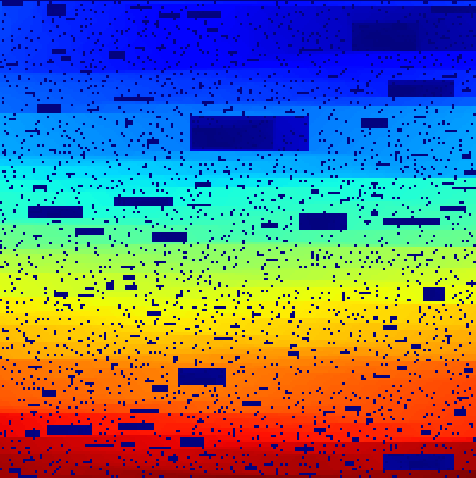}
	\caption{\emph{On the left:} a final configuration of the anisotropic model on $[40]^2$. Note that not all stable shapes are rectangles. \emph{On the right:} a final configuration on $[200]^2$ with $p=0.085$, where the color of each site represents the time it became infected. Blue sites became infected first, red sites last. }\label{fig:pictures}
\end{figure}

The anisotropic bootstrap percolation model is a specific two-dimensional process in the family described above. To be precise, set $d = 2$ and
\[
	\Nani \; := \; \big\{ (-2,0), (-1,0), (0,-1), (0,1), (1,0), (2,0) \big\}
\]
or graphically,\\
\begin{center}
\begin{tabular}{ccccccc}
	& & & & $\bullet$& & \\
	$\mathcal{N}_{\sss (1,2)}$ & $=$ & $\bullet$ & $\bullet$ & 0  & $\bullet$& $\bullet$\\
	& & & & $\bullet$& & \\
\end{tabular}
\end{center}
and set $r = 3$. $\Nani$ is sometimes called the ``$(1,2)$-neighbourhood'' of the origin. See Figure \ref{fig:pictures} for an illustration of the behaviour of the anisotropic model.

The main result of this paper is the following theorem:\footnote{
Throughout this paper we will use the standard Landau order notation: either for all $x$ sufficiently large or sufficiently small, depending on the context,
\begin{itemize}
	\item $f(x) = O(g(x))$ if there exists $C>0$ such that $|f(x)| \le C g(x)$,
	\item $f(x) = \Omega(g(x))$ if  there exists $c>0$ such that $|f(x)| \ge c g(x)$,
	\item $f(x) = \Theta(g(x))$ if $f(x) =O(g(x))$ and $f(x) = \Omega(g(x))$,
	\item $f(x) = o(g(x))$ if $| f(x)| /| g(x)| \to  0$.
\end{itemize}
}
\begin{thm}\label{thm:pc}
The critical probability of the anisotropic bootstrap percolation model satisfies
\begin{equation}\label{eq:thm:pc}
p_c\big( [L]^2, \Nani, 3 \big) \, = \, \frac{\log \log L}{12 \log L} \Big( \log \log L - 4 \log \log \log L + 2 \log\frac{9\e}{2} \pm o(1) \Big).
\end{equation}
\end{thm}

To put this theorem in context, let us recall some of the previous results obtained for bootstrap processes in two dimensions. The archetypal example of a bootstrap percolation model is the ``two-neighbour model'', that is, the process with neighbourhood 
\[
	\Nnn \, := \, \big\{ (-1,0), (0,-1), (0,1), (1,0) \big\}
\]
and $r = 2$. The strongest known bounds are due to Gravner, Holroyd, and Morris~\cite{GraHol08,GraHolMor12,Mor14}, who, building on work of Aizenman and Lebowitz~\cite{AizLeb88} and Holroyd~\cite{Hol03}, proved that 
\begin{equation}\label{e:2ndstandard}
	p_c\big( [L]^2, \Nnn, 2 \big) \, = \, \frac{\pi^2}{18 \log L} - \Theta\bigg( \frac{1}{(\log L)^{3/2}} \bigg).
\end{equation}

The anisotropic model was first studied by Gravner and Griffeath~\cite{GraGri96} in 1996.
In 2007, the second and third authors~\cite{EntHul07} determined the correct order of magnitude of $p_c$. More recently, the first and second authors~\cite{DumEnt13} proved that the anisotropic model exhibits a sharp threshold by determining the first term in~\eqref{eq:thm:pc}. 

The ``Duarte model'' is another anisotropic model that has been studied extensively \cite{Dua89,Mou95,BolDumMorSmi16}. The Duarte model has neighbourhood 
\[
	\Ndu \, = \, \big\{ (-1,0), (0,-1), (0,1) \big\}
\]
and $r = 2$. The sharpest known bounds here are due to the Bollob\'as, Morris, Smith, and the first author \cite{BolDumMorSmi16}:
\[
	p_c\big(\ZZ_L^2, \Ndu,2 \big) \, = \,  \frac{(\log \log L)^2}{8\log L}(1 \pm o(1)).
\]
Although the Duarte model has the same first order asymptotics for $p_c$ as the anisotropic model (up to the constant), the behaviour is very different. In particular, the Duarte model has a ``drift'' to the right: clusters grow only vertically and to the right. This asymmetry has severe consequences for the analysis of the model (especially for the shape of critical droplets).

The ``$r$-neighbour model'' in $d$ dimensions generalises the standard (two-neighbour) model described above. In this model, a vertex of $\ZZ^d$ is infected by the process as soon as it acquires at least $r$ already-infected nearest neighbours. Building on work of Aizenman and Lebowitz~\cite{AizLeb88}, Schonmann~\cite{Sch90}, Cerf and Cirillo~\cite{CerCir99}, Cerf and Manzo~\cite{CerMan02}, Holroyd~\cite{Hol03} and Balogh, Bollob\'as and Morris~\cite{BalBolMor09,BalBolMor10}, the following sharp threshold result for all non-trivial pairs $(d,r)$ was obtained by Balogh, Bollob\'as, Morris, and the first author~\cite{Baletal12}: for every $d \ge r \ge 2$, there exists an (explicit) constant $\lambda(d,r) > 0$ such that
\[
	p_c\big( [L]^d, \Nnnd, r \big) \, = \, \bigg( \frac{\lambda(d,r) \pm o(1)}{\log_{(r-1)} L} \bigg)^{d-r+1}.
\]
(Here, and throughout the paper, $\log_{\sss (k)}$ denotes a $k$-times iterated logarithm.)

Finally, we remark that much weaker bounds (differing by a large constant factor) have recently been obtained for an extremely general class of two-dimensional models by Bollob\'as, Morris, Smith, and the first author~\cite{Boletal14}, see Section~\ref{sec:history}, below. Moreover, stronger bounds (differing by a factor of $1 + o(1)$) were proved for a certain subclass of these models (including the two-neighbour model, but not the anisotropic model) by the first author and Holroyd~\cite{DumHol12}. 
 
 Although various other specific models have been studied (see e.g.~\cite{BriMah12,BriMahMel13,HolLigRom04}), in each case the authors fell very far short of determining the second term.
 
 \subsection{The bootstrap percolation paradox}
In \cite{Hol03} Holroyd for the first time determined sharp first order bounds on $p_c$ for the standard model, and observed that they were very far removed from numerical estimates: $\pi^2/18 \approx 0.55$, while the same constant was numerically determined to be $0.245 \pm 0.015$ on the basis of simulations of lattices up to $L = 28800$ \cite{AdlStaAha89}. This phenomenon became known in the literature as the \emph{bootstrap percolation paradox,} see e.g.\ \cite{Gra03,DeGetal04,GraHol08,AdlLev03}. 

An attempt to explain this phenomenon goes as follows: if the convergence of $p_c$ to its first-order asymptotic value is extremely slow, while for any fixed $L$ the transition around $p_c$ is very sharp, then it may appear that $p_c$ converges to a fixed value long before it actually does.

This indeed appears to be the case. The first rigorisation of the ``extremely slow convergence'' part of this argument appears in \cite{GraHol08}, for a model related to bootstrap percolation.
Theorem \ref{thm:pc} gives another unambiguous illustration of extremely slow convergence for a bootstrap percolation model: the second term in \eqref{eq:thm:pc} is actually larger than the first while 
\[
	4 \log \log \log L > \log \log L,
\]
which holds for all $L$ in the range $66 < L < 10^{2390}$. Moreover, the second term does not become negligible (smaller than $1\%$ of the first term, say) until $L > 10^{10^{1403}}$. On relatively small lattices, even the third term makes a significant contribution to $p_c$: it is larger than the first term when $L < 10^{60}$ and larger than the second term when $L < 10^{13}$.

The ``sharp transition'' part of the argument has also been made rigorous: for the standard model, an application of the Friedgut-Kalai sharp-threshold theorem \cite{BalBol03} tells us that the ``$\varepsilon$-window of the transition''\footnote{The $\varepsilon$-window denotes the difference between the value of $p_{\vep}$ where $[L]^2$ is internally filled with probability $\vep$, and the value $p_{1-\vep}$ where this probability equals $1-\vep$. In other words, the $\varepsilon$-window tells us how sharp the metastable transition is.} is 
\[
p_{1 - \varepsilon }([L]^2,\Nnn,2) -p_{\varepsilon}([L]^2,\Nnn,2) = O\left(\frac{\log \log L}{ {\log}^2 L}\right).
\]
So the $\varepsilon$-window is much smaller than the second order asymptotics in \eqref{e:2ndstandard}.

For the anisotropic model a similar analysis \cite{Boe11} yields that the $\varepsilon$-window satisfies
\[
p_{1- \varepsilon}([L]^2,\Nani,3)-p_{\varepsilon}([L]^2,\Nani,3) = O\left( \frac{\log^3 \log L}{ \log^2 L } \right),
\]
which is again much smaller than the second and third order asymptotics in Theorem \ref{thm:pc}.
So our analysis supports the above explanation of the bootstrap percolation paradox.

\subsection{Universality}\label{sec:history}

Recently, a very general family of bootstrap-type processes was introduced and studied by Bollob\'as, Smith and Uzzell~\cite{BolSmiUzz12a}. To define this family, let $\U = \{X_1,\ldots,X_m\}$ be a finite collection of finite subsets of $\ZZ^d \setminus \{0\}$, and define the corresponding bootstrap operator by setting
\[
	\B_{\sss \U}(\S) \, = \, \S \cup \big\{ v \in \ZZ^d \,:\, v + X \subset \S \;\textup{ for some } X \in \U \big\}
\]
for every set $\S \subset \ZZ^d$. It is not hard to see that all of the bootstrap processes described above can be encoded by such an `update family' $\U$, and in fact this definition is substantially more general. The key discovery of~\cite{BolSmiUzz12a} was that in two dimensions the class of such monotone cellular automata can be elegantly partitioned\footnote{This is the partitioning: We say a direction $u \in 
S^1$ is \emph{stable} if $\mathbb{H}_{u}$, the discrete half-plane that is orthogonal to $u$, satisfies $\langle \mathbb{H}_{u} \rangle = \mathbb{H}_{u}$. A family is \emph{supercritical} if there exists an open semicircle in $S^1$ containing no stable direction, and it is \emph{subcritical} if every open semicircle contains an infinite number of stable directions. It is \emph{critical} otherwise.} into three classes, each with completely different behaviour. More precisely, for every two-dimensional update family $\U$, one of the following holds:
\begin{itemize}
\item $\U$ is ``supercritical'' and has polynomial critical probability.\vspace{0.1cm}
\item $\U$ is ``critical'' and has poly-logarithmic critical probability.\vspace{0.1cm}
\item $\U$ is ``subcritical'' and has critical probability bounded away from zero. 
\end{itemize}
We emphasise that the first two statements were proved in~\cite{BolSmiUzz12a}, but the third was proved slightly later, by Balister, Bollob\'as, Przykucki and Smith~\cite{Baletal13}. Note that the critical class includes the two-neighbour, anisotropic and Duarte models (as well as many others, of course). For this class a much more precise result was recently obtained by Bollob\'as, Morris, Smith, and the first author~\cite{Boletal14}. In order to state this result, let us first (informally) define a two-dimensional update family to be ``balanced'' if its growth is asymptotically two-dimensional\footnote{In other words, in a balanced model the critical droplet is a polygon, all of whose sides have the same length up to a constant factor. For the precise definition, which is somewhat more technical, see~\cite{Boletal14}.} (like that of the two-neighbour model), and ``unbalanced'' if its growth is asymptotically one-dimensional (like that of the anisotropic and Duarte models). The following theorem was proved in~\cite{Boletal14}. 

\begin{thm}\label{thm:BDMS}
Let $\U$ be a critical two-dimensional bootstrap percolation update family. There exists $\alpha = \alpha(\U) \in \N$ such that the following holds:
\begin{enumerate}
\item[$(a)$] If $\U$ is balanced, then\footnote{Here $\ZZ_L^2$ denotes the discrete two-dimensional $L~\times~L$ torus, and $p_c\big( \ZZ_L^2,\U \big)$ is defined as in~\eqref{def:pc}. We consider the torus since in general undesirable complications may arise due to boundary effects or strongly asymmetrical growth.}
\[
	p_c\big( \ZZ_L^2,\U \big) = \Theta \bigg( \frac{1}{(\log L)^{1/\alpha}} \bigg).
\]
\item[$(b)$] If $\U$ is unbalanced, then
\[
	p_c\big( \ZZ_L^2,\U \big) = \Theta \bigg( \frac{(\log \log L)^2}{(\log L)^{1/\alpha}} \bigg).
\]
\end{enumerate}
\end{thm}
Theorem \ref{thm:BDMS} thus justifies our view of the anisotropic model as a canonical example of an unbalanced model.

\subsection{Internally filling a critical droplet}\label{sec:ISsub}

As usual in (critical) bootstrap percolation, the key step in the proof of Theorem~\ref{thm:pc} will be to obtain very precise bounds on the probability that a ``critical droplet'' $R$ is \emph{internally filled}\footnote{This notion is often referred to as ``internally spanned'' (especially in the older literature).} (IF), i.e., that $R \subset \< \S \cap R \>$. We will prove the following bounds: 
\begin{thm}\label{thm:IS}
Let $p > 0$ and $x,y \in \N$ be such that $1/p^2 \le x \le 1/p^5$ and $\frac{1}{3p} \log \frac{1}{p} \le y \le  \frac{1}{p} \log \frac{1}{p}$, and let $R$ be an $x \times y$ rectangle. Then 
\emph{
\[
\Pp \big( R \text{ is internally filled} \big) \, = \, \exp\left( -\frac{1}{6p} \left( \log \frac{1}{p} \right)^2 + \left( \frac{1}{3} \log\frac{8}{3\e} \pm o(1) \right) \frac{1}{p} \log \frac{1}{p} \right).
\]
}
\end{thm}

The alert reader may have noticed the following surprising fact: we obtain the first three terms of $p_c( [L]^2, \Nani, 3 )$ in Theorem~\ref{thm:pc}, despite only determining the first two terms of $\log \Pp ( R$ is IF$)$ in Theorem~\ref{thm:IS}. We will show how to formally deduce Theorem~\ref{thm:pc} from Theorem~\ref{thm:IS} in Section~\ref{sec:pc}, but let us begin by giving a brief outline of the argument.  

To slightly simplify the calculations, let us write 
\[
	C_1 : = \frac{1}{12} \qquad \text{ and }\qquad C_2 := \frac{1}{6} \log\frac{8}{3 \e}.
\]
We claim (and will later prove) that $p_c = p_c( [L]^2, \Nani, 3 )$ is essentially equal to the value of $p$ for which the expected number of internally filled critical droplets in $[L]^2$ is equal to $1$ (the idea being that a critical droplet with size as in Theorem \ref{thm:IS} will keep growing indefinitely with probability very close to one). We therefore have
\[
	L^2 \, \approx \, \exp\left( \frac{2C_1}{p_c} \left( \log \frac{1}{p_c} \right)^2 - \frac{2C_2}{p_c} \log \frac{1}{p_c} \right),
\]
and hence
\[
	p_c \, \approx \, \frac{C_1}{\log L} \left( \log \frac{1}{p_c} \right)^2 - \frac{C_2}{\log L} \log \frac{1}{p_c}.
\]
Iterating the right-hand side gives
\[
	p_c \, \approx \, \frac{C_1}{\log L} \left( \log \log L - 2 \log \log \frac{1}{p_c} - \log C_1 \right)^2 - \frac{C_2}{\log L} \left( \log \log L - 2 \log \log \frac{1}{p_c} \right).
\]
Upon using the approximation $\log\log(1/p_c) \approx \log\log\log L$ and multiplying out, this reduces to
\[
	p_c \, \approx \, \frac{C_1(\log \log L)^2}{\log L}  - \frac{4C_1\log \log L \log\log\log L}{\log L} - \big( C_2 + 2C_1\log C_1 \big) \frac{\log \log L}{\log L},
\]
which is what we hope to prove. Thus we obtain three terms for the price of two.

\subsection{A generalisation of the anisotropic model}\label{sec:onebee}

One natural way to generalise the anisotropic model is to consider, for each $b > a \ge 1$, the neighbourhood
\[
	\Ngen \; = \; \left(\big\{ (0,y) \in \ZZ^2 : -a \le y \le a \big\} \cup \big\{ (x,0) \in \ZZ^2 : -b \le x \le b \big\}\right) \setminus \{ (0,0) \}.
\]
It follows from Theorem~\ref{thm:BDMS} that 
\[
	p_c\big( \ZZ_L^2, \Ngen,r) \, = \, \Theta\bigg( \frac{(\log \log L)^2}{(\log L)^{1/\alpha}} \bigg),
\]
where $\alpha = r - b$, for each $b +1 \le r \le a + b$.\footnote{The value of $\alpha$ follows from \cite[Definition 1.2]{Boletal14}. Furthermore, if $r \le b$ then the model is supercritical, so $p_c\big( \ZZ_L^2, \Ngen,r) \le L^{-c}$ for some $c > 0$, and if $r > a + b$ then the model is subcritical, so $p_c\big( \ZZ_L^2, \Ngen,r) > c'$ for some $c' > 0$.} The arguments developed in \cite{DumEnt13} can be applied to prove that the leading order behaviour of $p_c$ for the $(1,b)$-model is\footnote{This is contrary to the claim in \cite[Section 1]{DumEnt13}. See also \cite{DumEnt15}, the erratum to \cite{DumEnt13}.}
\[
	p_c\big( [L]^2, \Ncal_{\sss (1,b)}, b+1 \big) \, = \, \bigg( \frac{(b-1)^2}{4(b+1)} \pm o(1) \bigg) \frac{(\log \log L)^2}{\log L}.
\]

Combining the techniques of~\cite{DumEnt13} with those introduced in this paper, it is possible to prove the following stronger bounds:

\begin{thm}\label{thm:onebee}
Given $b \ge 2$, set 
\[
	C(b) = \frac{2}{b-1} \log \left(\frac{{2b \choose b}  - \frac{2b-1}{b+1} {2b-2 \choose b} -1 + {2b \choose b-1} - {2b-2 \choose b-3}}{(b+1) \e}\right) + 2 \log \left(\frac{(b-1)^2}{4 (b+1)}\right).
\]
Then
\begin{equation}\label{eq:thm:pcb}
	p_c\big( [L]^2, \Ncal_{\sss (1,b)}, b+1 \big) \, = \, \frac{(b-1)^2}{4(b+1)} \frac{\log \log L}{\log L} \Big( \log \log L - 4 \log \log \log L - C(b) \pm o(1) \Big).
\end{equation}
\end{thm}

Note that in the case $b = 2$ this reduces to Theorem~\ref{thm:pc}. We remark that Theorem~\ref{thm:onebee} follows from a corresponding generalisation of Theorem~\ref{thm:IS}, with the constants $\frac16$ and $\frac13 \log \frac{8}{3\e}$ replaced by 
\[
	\frac{(b-1)^2}{2(b+1)} \qquad \text{ and } \qquad \frac{b-1}{b+1} \log \left( \frac{{2b \choose b}  - \frac{2b-1}{b+1} {2b-2 \choose b} -1 + {2b \choose b-1} - {2b-2 \choose b-3}}{(b+1)\e}\right),
\]
respectively.

\begin{figure}
	\includegraphics[width= 0.7\textwidth]{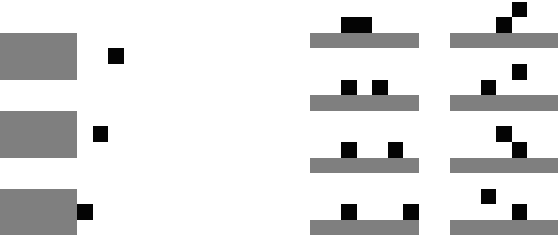}
	\caption{On the left: the three relative positions of an infected site that can cause horizontal growth of the grey rectangle. On the right: the 8 pairs of sites (up to shifts) that can cause vertical growth of the grey rectangle.}
	\label{fig:8pairs}
\end{figure}

We will not prove Theorem~\ref{thm:onebee}, since the proof is conceptually the same as that in the case $b = 2$, but requires several straightforward but lengthy calculations that might obscure the key ideas of the proof. It is, however, not too hard to see where the numerical factors come from:

A droplet grows horizontally in the $(1,b)$-model as long as it does not occur that the $b+1$ consecutive columns to its left and/or right do not contain an infected site. And it grows vertically as long as there are $b$ sites in a ``growth configuration'' somewhere above and/or below. There are 
\begin{equation}\label{e:1bnumber}
	{2b \choose b}  - \frac{2b-1}{b+1} {2b-2 \choose b} -1 + {2b \choose b-1} - {2b-2 \choose b-3}
\end{equation}
such configurations. Indeed, there are ${2b \choose b}$ different ways of finding $b$ infected sites inside $\Ncal_{\sss (1,b)} \setminus (\{(0,-1) ,(0,1)\})$. Of these, $\sum_{i=2}^b {2b-i \choose b} +1 = \frac{2b-1}{b+1} {2b-2 \choose b} +1$ are right-shifts of another configuration (e.g.\ for $\Ncal_{\sss (1,2)}$ the choices $\bullet \circ 0 \bullet \circ$ and $\circ \bullet 0 \circ \bullet$ count as a single growth configuration), so their contribution must be subtracted.
If $(0,1)$ is occupied, there are ${2b \choose b-1}$ ways of placing the other $b-1$ sites in $\Ncal_{\sss (1,b)} \setminus (\{(0,-1) ,(0,1)\})$. None of these are shift invariant, but some of them cannot grow to fill the entire row. Indeed, when $b \ge 3$, configurations where $(-b,0), (0,1),$ and $(b,0)$ are infected do not cause the row to fill up. Therefore, we must subtract ${2b-2 \choose b-3}$. This explains \eqref{e:1bnumber}. See Figure~\ref{fig:8pairs} for growth configurations of the case $b = 2$. 
Finally, it takes $b-1$ more infected sites for a rectangle to grow a row than it does to grow a column, which explains the remaining factors $b-1$ in~\eqref{eq:thm:pcb}.

\subsection{Comparison with simulations}
One might be tempted to hope that the third-order approximation of $p_c$ in Theorem~\ref{thm:pc} is reasonably good already for lattices that a computer might be able to handle. Simulations indicate that this is not the case. Indeed, for lattices with $L \le 10,\!000,$ the third-order approximation is even farther from to the simulated values than the first-order approximation (and recall that the second-order approximation is negative here). We believe that this should not be surprising, because it is not at all obvious that the fourth order term should be significantly smaller: careful inspection of our proof suggests that the $o(\frac1p \log \frac1p)$ term in Theorem~\ref{thm:IS} is at most $O(\frac1p \log \log \frac1p)$. Although we do not prove this, we have no reason to believe that a correction term of that order does not exist. Even if we suppose that the third order correction in Theorem~\ref{thm:IS} can be sharply bounded by $C_3 /p$, say, so that we would have the bound
\[
	\Pp \big( R \text{ is internally filled} \big) \, \stackrel{?}{=} \, \exp\left( -\frac{2C_1}{p} \left( \log \frac{1}{p} \right)^2 +  \frac{2 C_2}{p} \log \frac{1}{p} + \frac{C_3 \pm o(1)}{p} \right),
\]
for critical droplets instead, then a computation like the one in Section~\ref{sec:ISsub} above suggests that this would yield
\begin{multline*}
	p_c([L]^2,\Nani,3) \stackrel{?}{=} \frac{\left(\log \log L\right)^2}{12\log L} \, - \, \frac{(\log \log L )( \log \log \log L)}{ 3\log L}
	 +  \frac{\left(\log \frac{9}{2}+1\right) \log \log L}{6\log L} \\
	 + \frac{(\log \log \log L)^2}{3 \log L} - \frac{(\log \frac92 +1) \log\log \log L}{3\log L}\\
	  + \frac{C_3 + \log \frac{1}{12} \left(2 + \log \frac{27}{16}\right) \pm o(1)}{12 \log L},
\end{multline*}
so the fourth, fifth, and sixth order terms of $p_c$ would also be comparable to the first for moderately sized lattices. Moreover, because of the extremely slow decay of these correction terms (e.g.\ $(\log \log 10^{10})^2 \approx 10$), it might be too optimistic to expect that one would be able to determine $C_3$ by fitting to the simulated values of $p_c$, if indeed $C_3$ exists.

\subsection{Comparison with the two-neighbor model}
Comparing Theorem~\ref{thm:pc} with the analogous result for the two-neighbor model, \eqref{e:2ndstandard}, it may seem remarkable how much sharper the former is than the latter. We believe the following heuristic discussion goes a way towards explaining this difference.

Both approximations of $p_c$ are proved using essentially the same critical droplet heuristic described above. Once a critical droplet has formed, the entire lattice will easily fill up. But filling a droplet-sized area is exponentially unlikely: it is essentially a large deviations event. The theory of large deviations tells us that if a rare event occurs, it will occur in the most probable way that it can. For filling a droplet, this means that one should find an optimal ``growth trajectory'': a sequence of dimensions from which a very small infected area (a ``seed'') steadily grows to fill up the entire droplet. For the anisotropic model, in \cite{DumEnt13}, the first and second authors determined this trajectory to be close to $x = \frac{\e^{3py}}{3p}$, where $x$ and $y$ denote the horizontal and vertical dimensions of the seed as it grows. This approximation was enough to yield the first term of $p_c$. In the current paper we establish tighter bounds of optimal trajectory around $x = \frac{\e^{3py}}{3p}$, allowing us to give the sharper estimate for the probability of filling a droplet in Theorem~\ref{thm:IS}. As we showed in Section~\ref{sec:ISsub} above, this correction is enough to obtain the first three terms of $p_c$ for the anisotropic model.

For the two-neighbor model, however, finding this optimal growth trajectory is not at all the challenge: by symmetry it is trivially $x=y$. The correction to $p_c$ that Gravner, Holroyd, and Morris determined in \cite{GraHol08,GraHolMor12,Mor14}, is instead due to the much smaller entropic effect of random fluctuations around this trajectory (see also the introduction of \cite{GraHol09} for a more detailed explanation of this effect). We believe that such fluctuations also influence $p_c$ for the anistropic model, but that their effect will be much smaller than the improvements that can still be made in controlling the precise shape of the optimal growth trajectory.

\subsection{About the proofs}
The proof of Theorem \ref{thm:pc} uses a rigorisation of the iterative determination of $p_c$ in Section \ref{sec:ISsub} above, combined with Theorem \ref{thm:IS} and the classical argument of Aizenman and Lebowitz \cite{AizLeb88}. 

The lower bound in Theorem \ref{thm:IS} is a refinement of the computation in \cite{DumEnt13}.

Most of the work of this paper goes into the proof of the upper bound of Theorem~\ref{thm:IS}. Like many recent entries in the bootstrap percolation literature, our proof centers around the ``hierarchies'' argument of Holroyd \cite{Hol03}. In particular, we sharpen the argument of \cite{DumEnt13} by incorporating the idea of ``good'' and ``bad'' hierarchies from \cite{GraHolMor12}, and by using very precise bounds on horizontal and vertical growth of infected rectangular regions.

The main new contributions of this paper (besides the iterative determination of $p_c$) can be found in Sections \ref{sec:kvert} and \ref{sec:Wbd}. 

In Section \ref{sec:kvert}, we introduce the notion of \emph{spanning time} (Definition~\ref{def:tau}), which characterises to a large extent the structure of configurations of vertical growth. We show that if the spanning time is $0$, then such structures have a simple description in terms of paths of infected sites, whereas if the spanning time is not $0$, then this description can still be given in terms of paths, but these paths now also involve more complex arrangements of infected sites. We call such arrangements \emph{infectors} (Definition~\ref{def:infector}), and show that they are sufficiently rare that their contribution does not dominate the probability of vertical growth.

In Section \ref{sec:Wbd} we generalise the variational principle of Holroyd \cite{Hol03} to a more general class of growth trajectories. This part of the proof is intended to be more widely applicable than the current anisotropic case, and is set up to allow for precise estimates.

\subsection{Notation and definitions} 
A \emph{rectangle} $[a,b] \times [c,d]$ is the set of sites in $\ZZ^2$ contained in the Euclidean rectangle $[a,b] \times [c,d]$. 
For a finite set $\Q \subset \ZZ^2$, we denote its \emph{dimensions} by $(\xb(\Q), \yb(\Q))$, where $\xb(\Q) = \max\{a_1 - b_1 +1 \, : \, \{(a_1,a_2), (b_1, b_2)\} \in \Q \times \Q\}$, and similarly, $\yb(\Q) = \max\{a_2 - b_2 +1\, : \, \{(a_1,a_2), (b_1, b_2)\} \in \Q \times \Q\}$. So in particular, a rectangle $R = [a,b] \times [c,d]$ has dimensions $(\xb(R), \yb(R)) = (|[a,b] \cap \ZZ|, |[c,d] \cap \ZZ|)$. Oftentimes, the quantities that we calculate will only depend on the size of $R$, and be invariant with respect to the position of $R$. In such cases, when there is no possible confusion, we will write $R$ with $\xb(R)=x$ and $\yb(R)=y$ as $[x] \times [y]$. A \emph{row} of $R$ is a set $\{(m,n) \in R \,:\, n=n_0\}$ for some fixed $n_0$. A \emph{column} is similarly defined as a set $\{(m,n) \in R \,:\, m =m_0\}$.  We sometimes write $[a,b] \times \{c\}$ for the row $\{ (m,c) \in \ZZ^2 \, :\, m \in [a,b] \cap \ZZ\}$, and use similar notation for columns.

We say that a rectangle $R = [a,b] \times [c,d]$ is \emph{horizontally traversable} (hor-trav) by a configuration $\S$ if 
\[
	R \subset \langle (R \cap \S) \cup ([a-2,a-1] \times [c,d])\rangle.
\]
That is, $R$ is horizontally traversable if the rectangle becomes infected when the two columns to its left are completely infected.
Under $\Pr_p$, this event is equiprobable to the event that $R \subset \langle (R \cap \S) \cup ([b+1, b+2] \times [c,d])\rangle$, and more importantly, it is equivalent to the event that $R$ does not contain three or more consecutive columns without any infected sites and the rightmost column contains an infected site.

We say that $R$ is \emph{up-traversable} (up-trav) by $\S$ if
\[
	R \subset \langle (R \cap \S) \cup ([a,b] \times \{c-1\})\rangle.
\]
That is, $R$ becomes entirely infected when all sites in the row directly below $R$ are infected. Similarly, we say that $R$ is \emph{down-traversable} by $\S$ if $R \subset \langle (R \cap \S) \cup ([a,b] \times \{d+1\})\rangle$. Again, under $\Pr_p$ up and down traversability are equiprobable, so we will only discuss up-traversability. If $\S$ is a random site percolation, then we simply say that $R$ is horizontally- or up- or down-traversable.

Given rectangles $R \subset R'$ we write $\left\{R \growto R' \right\}$ for the event that the dynamics restricted to $R'$ eventually infect all sites of $R'$ if all sites in $R$ are infected, i.e., for the event that $R' = \langle (\S \cap R') \cup R\rangle$. 

We will frequently make use of two standard correlation inequalities: The first is the \emph{Fortuin-Kasteleyn-Ginibre inequality} (FKG-inequality), which states that for increasing events $A$ and $B$, $\Pp(A \cap B) \ge \Pp(A) \Pp(B)$. The second is the \emph{van den Berg-Kesten inequality} (BK-inequality), which states that for increasing events $A$ and $B$, $\Pp(A \circ B) \le \Pp(A) \Pp(B)$, where $A \circ B$ means that $A$ and $B$ occur disjointly (see \cite[Chapter 2]{Gri99} for a more in-depth discussion).

\subsection{The structure of this paper} In Section~\ref{sec:lower} we state two key bounds, Lemmas~\ref{21} and \ref{22}, giving primarily lower bounds on the probabilities of horizontal and vertical growth of an infected rectangular region, and we use them to prove the lower bound of Theorem~\ref{thm:IS}. In Section~\ref{sec:kvert} we prove a complementary upper bound on the vertical growth of infected rectangles, Lemma~\ref{lem:kvc}. In Section~\ref{sec:lemWproof} we prove Lemma~\ref{lem:W}, which combines the upper bounds on horizontal and vertical growth from Lemmas~\ref{21} and \ref{lem:kvc}. This lemma is crucial for the upper bound of Theorem \ref{thm:IS}. We prove the upper bound of Theorem~\ref{thm:IS} in Section~\ref{sec:upper}, subject to a variational principle, Lemma~\ref{lem:Wbd}, that we prove in Section \ref{sec:Wbd}. 
Finally, in Section~\ref{sec:pc} we use Theorem~\ref{thm:IS} to prove Theorem~\ref{thm:pc}. 

\subsection*{Acknowledgments} The authors would like to thank Robert Morris for his involvement in the earlier stages of the project, and for the many crucial insights he provided. We thank the anonymous referee for their careful reading and comments. The third author would like to thank Robert Fitzner for useful discussions about computer simulations.

The first author was supported by the IDEX Chair funded by Paris-Saclay and the NCCR SwissMap funded by the Swiss NSF.

The third author was supported by the Netherlands Organisation for Scientific Research 
(NWO) through Gravitation-grant {\sc networks}-024.002.003.

\section{The lower bound of Theorem \ref{thm:IS}}\label{sec:lower}
Recall that $C_1 = \frac{1}{12}$ and $C_2 = \frac{1}{6} \log \frac{8}{3 \e}$.

\begin{prop}\label{upper}
Let $p>0$ and $\frac{1}{p^2} \le x \le \frac{1}{p^5}$ and $\frac{1}{3p} \log \frac1p \le y \le \frac{1}{p^5}$. Then \emph{
\[
	\Pp ([x] \times [y] \text{ is IF}) \; \ge \; \exp\left(-\frac{2C_1}{p} \log^2 \frac 1p + \left( 2C_2 - o(1) \right)\frac1p \log \frac1p \right).
\]}
\end{prop}
Note that the upper bound on $y$ is different from the bound in Theorem~\ref{thm:IS}.

For the proof it suffices to show that there exists a subset of configurations that has the desired probability. We choose a subset of configurations that follow a typical ``growth trajectory'': configurations that contain a small area that is locally densely infected (a \emph{seed}). We bound the probability that such a seed will grow a bit (which is likely), and then a lot more (which is exponentially unlikely), until the infected region reaches a size where the growth is again very likely, because the boundary of the infected region is large and the dynamics depend only on the existence of infected sites on the boundary, not on their number.

To prove this proposition we will need bounds on the probability that a rectangle becomes infected in the presence of a large infected cluster on its boundary. We state two lemmas that achieve this, which are improvements upon \cite[Lemmas 2.1 and 2.2]{DumEnt13}.

\begin{lemma}\label{21}
For any rectangle $[x] \times [y]$,
\[
	\e^{ -x f(p,y)} \; \le \; \Pp \left( [x] \times [y] \textup{ is hor-trav} \right) \; \le \; \e^{-(x-2)f(p,y)},
\]
where $f(p,y) := -\log( \alpha (1-(1-p)^y))$ and where $\alpha(u)$ is the positive root of the polynomial 
\begin{equation}\label{e:poly}
	X^3-uX^2-u(1-u)X-u(1-u)^2.
\end{equation}
Moreover, $f(p,y)$ satisfies the following bounds:
\begin{enumerate}
	\item when $p\rightarrow 0$ and $py\rightarrow \infty$,
	 \[
	 	f(p,y)=\e^{-3py}+\Theta(\e^{-4py}),
	\]
	\item when $y \ge \frac{2}{p} \log\log \frac1p$,
	\[
		f(p,y)=\e^{-3py}\left(1+\Theta\left(\log^{-2} (1/p)\right)\right),
	\]
	\item when $p \to 0$, $y \to \infty$, and $(1-p)^y \to 1$,
	\[
		f(p,y) \ge \tfrac12 p y - 3 p^2 y^2.
	\]
\end{enumerate}
\end{lemma} 
\proof
From \cite[Lemma 2.1]{DumEnt13}\footnote{Note that in the proof of \cite[Lemma 2.1]{DumEnt13} there are a number of (unimportant) sign errors.} we know that 
\[
	\alpha \left(1-(1-p)^y \right)^{x} \;\le\; \Pp \left([x] \times [y] \textup{ is hor-trav} \right)\; \le \; \alpha \left(1-(1-p)^y \right)^{x-2}.
\]
When $u$ is close to $1$, $X = \e^{-(1-u)^3}$ is an approximate solution for the positive root, since
\[
	\e^{-3(1-u)^3} -u\e^{-2(1-u)^3} -u(1-u) \e^{-(1-u)^3} -u(1-u)^2 = \Theta((1-u)^4 ).
\]
So, as $p \to 0$ and $py \to \infty$,
\[
	-\log \alpha(1-(1-p)^y) = (1-p)^{3y} + \Theta((1-p)^{4y}) = \e^{-3py} + \Theta(\e^{-4py}).
\]
This establishes (a) and (b) simply follows.
\medskip

To prove (c), recall Rouch\'e's Theorem (see e.g.\ \cite[Theorem 10.43]{Rud87}), which states that if two functions $g(z)$ and $h(z)$ are holomorphic on a bounded region $U \subset \mathbb{C}$ with continuous boundary $\partial U$ and satisfy $|g(z) - h(z)| < |g(z)|$ for all $z \in \partial U$, then $g$ and $h$ have an equal number of roots on $U$. Applying Rouch\'e's Theorem with $h(z) = a_0 + a_1 z + a_2 z^2 + a_3 z^3$ and $g(z) = a_0$,  it follows that the moduli of the roots of $h(z)$ are all bounded from below by $|a_0| /(|a_0| + \max\{|a_1|, |a_2|, |a_3|\})$. Applying this bound to \eqref{e:poly} we find that when $u > 0$ is sufficiently small,
\[
	\alpha(u) \ge \frac{u(1-u)^2}{u(1-u)^2+1} \ge u -3 u^2,
\]
where the second inequality is due to a series expansion around $u=0$. (We remark that an explicit computation gives $\alpha(u) \ge u -3u^2$ for all $u >0$, but without relying on a computer this may take several pages to verify.)
Since we assumed $(1-p)^y \to 1$ we thus have
\[
	f(p,y) \ge (1-(1-p)^y) - 3(1-(1-p)^y)^2 \ge \tfrac12 py - 3 p^2 y^2,
\]
where we used $\frac12 py \le 1-(1-p)^y \le py$ for $p$ sufficiently small.
\qed

\begin{lemma}\label{22}
(a) If $p^2x$ is sufficiently small, then we have, for any rectangle $[x] \times [y]$,
\[
	\Pp \left([x] \times [y] \textup{ is up-trav} \right) \; \ge \; \exp\Big( y \log(8p^2x)\big(1+O(p^2x + p)\big) \Big).
\]

(b) As long as $ \frac{8 p^2 x}{5} \leq 1$ we have 
\[
	\Pp([x] \times [y] \textup{ is up-trav}) \; \ge \; \left(\frac{8 p^2 x}{5 \e} \right)^y.
\]
\end{lemma}

\proof
 We say that a rectangle is \emph{North-traversable} (N-trav) if the intersection of every row with $R$ contains a site $(a,b)$ such that $((a,b) + \Nani)\setminus \{(a, b-1)\}$ contains at least two infected sites. Observe that North-traversability implies up-traversability, so
 \[
 	\Pp([x] \times [y] \text{ is up-trav}) \ge \Pp([x] \times [y] \textup{ is N-trav}).
\]
We can similarly define South-traversability by requiring that the intersection of every row with $R$ contains a site $(a,b)$ such that $((a,b) + \Nani)\setminus \{(a, b+1)\}$ contains at least two infected sites. South-traversability implies down-traversability. Again, from a probabilistic point of view North- and South-traversability are equivalent, so we will henceforth only discuss North-traversability.

If $[x] \times[y]$ is North-traversable then for each of the $y$ rows there must exist an infected pair of sites $u$ and $v$ and a site $z$ in the row such that $u,v \in z+ \Nani$. 
By the FKG inequality we thus have the lower bound
\[
	\Pp([x] \times [y] \textup{ is N-trav}) \ge \Pp(\exists \text{ an infected pair for a row of length $x$})^y.
\]

For the proof of (a) we apply Janson's inequality \cite{Jan90}. The expected number of infected pairs immediately above an infected rectangle of width $x$ is at least $\mu~=~(8x - 16)p^2$. To see this, consider that up to translations there are $8$ possible pairs of infected sites above the rectangle that can infect the whole row, see Figure \ref{fig:8pairs} above. The variance is\footnote{For two positive sequences $a_n$ and $b_n$ we write $a_n \gg b_n$ when $a_n/b_n \to \infty$ and $a_n \ll b_n$ when $a_n/b_n \to 0$.} $\Delta = O(p^3x) \ll \mu$, so the probability that some pair is infected is at least 
\[
	1 - \exp\left( -\mu + \Delta/2 \right) \; \ge \; \left( 8p^2 x-O(p^3x + p^4x^2) \right),
\]
using the inequality $1 - \e^{-x} \ge x - x^2$ for $x \ge 0$.
\medskip

For the proof of (b) we use a cruder approximation: 
For $(a,b) \in [x] \times [y]$ let $A_{(a,b)}$ be the event that $(a,b)$ is the leftmost site of an infected pair as in Figure \ref{fig:8pairs}.
These pairs all have width at most $5$, so the probability that a row of length $x$ does not have an infected pair can be bounded from above by
\[
	(1-8p^2)^{\lfloor x /5 \rfloor} \le \exp\left(-\frac{8 p^2 x}{5}\right) \le 1-\frac{8 p^2 x}{5 \e}
\]
when $ \frac{8 p^2 x}{5} \le 1$. The claim follows. \qed
\medskip

\proof[Proof of Proposition \ref{upper}]
We start by constructing a seed. 
Let $r := \lfloor \frac2p \log\log \frac1p \rfloor$ and infect sites $(1,2i)$ and $(2,2i+1)$ for $2i\leq r$. The probability that a rectangle $[2] \times [r]$ is a seed is $p^r$.
Note that the infected sites internally fill $[2]\times[r]$.

The growth of the seed to a rectangle of arbitrary size can be divided into three stages:
\medskip

\emph{Stage 1.} By Lemma \ref{21}(a) the probability of finding a seed of size $r$ that will grow to size $\left[\e^{3 rp}/(3p)\right] \times [r]$ is about the same as the probability of just finding the seed, i.e.,
\begin{equation}\label{eq-ubseed}
	p^{r} \cdot \exp\left( - \frac{\e^{3rp}}{3p}\cdot \left(\e^{-3 r p}+O\left(\e^{-4 r p}\right)\right) \right)
	 \ge p^r \e^{-O(1/p)}.
\end{equation}
\medskip

\emph{Stage 2.} Next we bound the probability that the infected rectangle grows to size
\[
	R := \left[ \frac{1}{3 p^2} \right]\times\left[\frac1{3p} \log \frac 1p \right],
\]
that is, we want to bound
\begin{equation}\label{e:startstage2}
	\Pp \left(\left [\frac{\e^{3rp}}{3p}\right] \times \left[ r \right] \growto \left[ \frac{\e^{3mp}}{3p}\right] \times \left[m \right] \right),
\end{equation}
where $m:= \frac{1}{3p} \log \frac1p$. This is the bottleneck for the growth dynamics.
We bound \eqref{e:startstage2} by considering the growth in many small steps. In each such step, the rectangle will either infect an entire row above or below it, or it will infect an entire row to the left or right of it (with the help of infected sites on the boundary of the rectangle). Because vertical growth is less probable than horizontal growth, we will consider sequences where the rectangle grows by one vertical step, from height $\ell$ to $\ell+1$, followed by horizontal growth that infects many columns successively, with the rectangle growing from width $x_\ell$ to $x_{\ell+1}$ where $ x_\ell: = \frac{\e^{3 \ell p}}{3p}$.
That this choice is close to optimal can be seen in Section \ref{sec:Wbd} below, where a variational principle for the upper bound of Theorem~\ref{thm:IS} is derived.

Having divided the growth into steps, we can bound \eqref{e:startstage2} from below using the FKG-inequality:
\begin{equation}\begin{split}\label{e:3prod}
	\Pp \Biggl(\left [\frac{\e^{3rp}}{3p}\right] \times \left[ r \right] \growto  \left[ \frac{\e^{3mp}}{3p}\right] \times \left[m \right] \Biggr) 
	\ge& \prod_{\ell=r}^{m} \Pp \left(\left[x_\ell\right] \times \left[ \ell \right] \growto \left[ x_{\ell+1}\right] \times \left[\ell\right]\right)\\ 
	&\times \prod_{\ell=r}^{m-r} \Pp \left( \left[x_\ell \right] \times \left[ \ell \right] \growto \left[ x_\ell \right] \times \left[\ell+1\right]\right)\\
	&\times \prod_{\ell=m-r+1}^{m} \Pp \left( \left[x_\ell \right] \times \left[ \ell \right] \growto \left[ x_{\ell}\right] \times \left[\ell+1\right]\right).
\end{split}\end{equation}
We bound these three products separately.	

It follows from Lemma \ref{21}(a) that the horizontal growth from width $x_\ell$ to $x_{\ell+1}$ occurs with probability approximately $1/\e$, i.e.,
\begin{equation}\label{e:4prod}\begin{split}
	\Pp \big( [x_\ell] \times [\ell] & \growto [x_{\ell+1}] \times [\ell]\big) \\
	& \ge \exp\Big(-\frac{1}{3p} \big( \e^{ 3  (\ell+1)p } - \e^{3  \ell p } \big)  \e^{- 3  \ell p }\big(1+O\big(\log^{-4/3} (1/p)\big)\big)\Big) \\
	& \ge \e^{-1 - o(1)}.
	\end{split}
\end{equation}
Therefore,
\begin{equation}\label{e:prod1}
		\prod_{\ell=r}^{m} \Pp \big( [x_\ell] \times [\ell] \growto [x_{\ell+1}] \times [\ell]\big) \ge \e^{-m (1+o(1))}.
\end{equation}

When $\ell \le m-r$, then $p^2 x_\ell \le \log^{-2} \frac 1p$, so we can apply Lemma \ref{22}(a) to bound
\[
	\Pp ([x_\ell] \times [\ell] \growto [x_\ell] \times [\ell+1]) \ge 8 p^2 x_\ell \, \e^{O\left(\log^{-4/3} \frac 1p\right)}.
\]
Therefore we can bound the second product in \eqref{e:3prod} from below by
\begin{equation} \label{e:prod2}
	\begin{split}
		\prod_{\ell=r}^{m-r} 8 p^2 x_\ell \e^{O\left(\log^{-4/3} \frac 1p\right)} & \ge \left(\frac{8 p}{3}\right)^{m-2r} \exp\left( 3p \sum_{\ell=r}^{m-r} \ell \right) \e^{(m-2r)O\left(\log^{-4/3} \frac 1p\right)}\\
			& =  \left(\frac{8 p}{3}\right)^{m-2r} \exp \left( \frac{3p}{2} \big((m-r) (m-r+1) - (r-1)r) \big) \right) \e^{o\left(m \right)}\\
			& =  p^{-r} \left(\frac{8 p}{3}\right)^{m-r} \exp \left( \frac{3p}{2} \left(m^2 - 2mr + m\right) \right) \e^{o\left(m\right)}\\
			& = p^{-r} \left(\frac{8 p }{3} \right)^{m-r} \exp \left(\frac{3 p}{2} (m^2 - 2mr)\right) \e^{o(m)}.
	\end{split}
\end{equation}

Using Lemma \ref{22}(b) we can similarly bound the third product from below by
\begin{equation}\label{e:prod3}
 	\begin{split}
		\prod_{\ell=m-r+1}^{m} \frac{8 p^2 x_\ell}{5 \e} &= \left(\frac{8p}{15 \e}\right)^r \exp\left(\sum_{\ell = m-r+1}^m 3 \ell p \right)\\
		&\ge \left(\frac{8 p}{3}\right)^{r} \exp \left( \frac{3p}{2} \big(m(m+1) - (m-r)(m-r+1) \big)\right)\left(\frac{1}{5\e}\right)^{r}\\
		& = \left(\frac{8p}{3}\right)^{r} \exp \left(\frac{3p}{2} 2mr\right) \e^{o(m)}.
	\end{split}
\end{equation}
Multiplying the bounds \eqref{e:prod1}, \eqref{e:prod2}, and \eqref{e:prod3}, and using that $m = \frac{1}{3p} \log \frac1p$, we get
\begin{equation}\label{e:stage2}
\begin{split}
\Pp \left(\left [\frac{\e^{3rp}}{3p}\right] \times \left[ r \right] \growto \left[ \frac{\e^{3mp}}{3p}\right] \times \left[m \right] \right) &\ge  p^{-r} \left(\frac{8 p}{3}\right)^m \exp\left(\frac{3 p}{2} m^2 -m\right)\e^{o(m)}\\
	& = p^{-r} \exp\left(\frac{3 p}{2} m^2 - m \log \frac1p + m \log \frac{8}{3} -m \right) \e^{o(m)}\\
	&= p^{-r} \exp\left(- \frac{1}{6p} \log^2 \frac1p + (1- o(1))\frac{1}{3p} \log \frac{8}{3\e}  \log \frac1p  \right).
\end{split}
\end{equation}
	
\emph{Stage 3.} The infected region can grow from $[\frac{1}{3p^2}] \times [m]$ to arbitrary size with good probability. 
Indeed, we claim that 
\begin{equation}\label{e:stage3}
	\Pp \left(\left[\frac{1}{3p^2} \right] \times [m] \growto R\right) \ge \e^{-O\left( 1/p \right)}.
\end{equation}
This bound is proved in \cite[proof of Proposition 2.4]{DumEnt13}.
We do not repeat the proof here, but let us indicate how this bound is established: Consider the case where the cluster first grows horizontally to width $1/p^2$. By Lemma \ref{21}(b) we have
\[
	\Pp \left(\left[\frac{1}{3p^2} \right] \times [m] \growto \left[\frac{1}{p^2}\right] \times [m] \right) \ge \exp\left( - \frac{2}{3 p^2} \cdot p (1 +o(1))\right)  =  \e^{- O\left(1/p\right)}.
\]
Now consider the case where it grows vertically, this time to height $3m$. This also occurs with probability at least $\e^{-O(1/p)}$. As the infected region gets larger, the probability that it keeps growing converges to $1$. The result is that \eqref{e:stage3} holds for any rectangle $R$ that is large enough, as long as the dimensions of $R$ are sufficiently balanced (which is guaranteed by the assumptions on $x$ and $y$).

Now, by the FKG-inequality, we can multiply the bounds from the three stages (i.e., \eqref{eq-ubseed}, \eqref{e:stage2}, and \eqref{e:stage3}) to complete the proof of Proposition~\ref{upper}.\qed
\medskip


\section{An upper bound on the probability of up-traversability}\label{sec:kvert}

The following bound is crucial for the proof of the upper bound of Theorem~\ref{thm:IS}. 
Recall from~\eqref{def:operator} the definition of the bootstrap operator $\B$, and recall that $\B^{(t)}(\S)$ is the $t$-th iterate of $\B$ with initial set $\S$, and that $\<\S\> = \lim_{t \to \infty} \B^{(t)}(\S)$. Recall that a rectangle $R = [1,x] \times [1,y]$ is said to be up-traversable by a set $\S$ if $R \subset \< (\S \cap R) \cup ([1,x] \times \{0\}) \>$, and that we write $\Pp$ to indicate that the elements of $\S$ are chosen independently at random with probability $p$.

\begin{lemma}\label{lem:kvc} 
Let $1 \le k \ll p^{-2/5}$ and let $R$ be a rectangle with dimensions $(x,y)$ such that $y < x$.
Then, for $p$ sufficiently small,
\[
	\Pp\big( R \text{ is up-traversable} \big) \le  \begin{cases} p^{-k} \e^{y/k} (24 p k^2 + 8p)^y & \text{ if } \quad  x < \frac{3 k^2}{p},\\
	p^{-k} \e^{y/k} \big( 8p^2 x + 8p \big)^{y} & \text{ if } \quad \frac{3 k^2}{p} \le x \le \frac{1}{p^2}.
	\end{cases}
\]
\end{lemma}

We will apply this lemma with $\frac1p \ll y \ll \frac{1}{p} \log^6 \frac1p \le  x$ and $k = \log^2 \frac1p$. Note that in this case the upper bound given by the lemma is not much larger than the lower bound given by Lemma~\ref{22}. In particular, for these choices of $x, y$ and $k$, the bound given by the lemma is of the form $\big( ( 8 + o(1)) p^2 x \big)^y$.

We begin the proof of Lemma~\ref{lem:kvc} with the following simple but important definition: let us say that a pair of sites $\P$ is a \emph{spanning pair} for the row $[a,b] \times \{c\}$ if 
\begin{equation}\label{def:spanning:pairs}
[a,b] \times \{c\} \subset \langle \P \cup [a,b] \times \{c-1\}\rangle.
\end{equation}
That is, $\P$ is a spanning pair for $[a,b] \times \{c\}$ if the row becomes infected when $\P$ and the row below it are infected.
Note that for each spanning pair $\P = \{u,v\}$ there exists $z \in [a,b] \times \{c\}$ such that $u,v \in z + (\Nani \setminus \{(0,-1)\})$, and thus that any spanning pair is a translate of one of the eight pairs on the right-hand side of Figure \ref{fig:8pairs}.

\begin{lemma}\label{lem:equivuptrav} 
Let $R$ be a rectangle such that $R$ has $\xb(R) \ge 2$ and $\yb(R) \ge 1$, and let $\S \subset R$. Then $R$ is up-traversable by~$\S$ if and only if $\< \S \>$ contains a spanning pair for every row of~$R$.
\end{lemma}

\proof
Suppose that $R = [a,b] \times [c,d]$ with $b-a \ge 1$ and $d-c \ge 0$.
It is easy to see that if $\< \S \>$ contains a spanning pair for every row of~$R$, then $R$ is up-traversable by $\S$: if $\<\S\>$ contains a spanning pair for the bottom row of $R$, then the whole row becomes infected, i.e., $[a,b] \times \{c\} \subset \< \S \cup [a,b] \times \{c\}\>$. And given that the bottom row is infected, the row above the bottom row must also become infected, since $\<\S\>$ also contains a spanning pair for it, i.e., $[a,b] \times \{c+1\} \subset \<\S \cup [a,b] \times \{c\}\>$. This argument can be repeated for all rows. 

It will therefore suffice to prove that the converse holds.
To do that, let $j \in [c,d]$ be the smallest $j$ such that $\< \S \>$ does not contain a spanning pair for the row $[a,b] \times \{j\}$. We claim that the set
\[
	 \big( \< \S \cup [a,b] \times \{j-1\} \> \setminus \< \S \> \big) \cap ([a,b] \times \{j\})
\] 
is empty. Indeed, suppose that for some $t \ge 1$ there exists a site $v$ such that 
\[
	v \in \Bcal^{(t)}\big(\S \cup ([a,b] \times \{j-1\})\big) \cap ([a,b] \times \{j\}),
\]
Then there must be a pair of already-infected sites in $\Nani(v) \cap ([a,b] \times [j,j+1])$ at time $t-1$. But this pair lies in $\< \S \>$, and thus is a spanning pair for the row $[a,b] \times \{j\}$, a contradiction. Now, since $\< \S \>$ does not contain a spanning pair for $[a,b] \times \{j\}$, this implies that $R \nsubseteq \< \S \cup ([a,b] \times \{c-1\}) \> $, as required.\qed
\medskip

We now make another important definition.

\begin{defn}\label{def:tau}
For each rectangle $R$  such that $\xb(R) \ge 2$ and $\yb(R) \ge 1$, and each set $\S \subset R$ such that $R$ is up-traversable by~$\S$, let $\A(\S) \subset \S$ be a minimum-size subset such that $R$ is up-traversable by $\A(\S)$. (If more than one such subset exists, then choose one according to some arbitrary rule.) Define the \emph{spanning time}
\[
	\tau \, = \, \tau(R, \S) \, := \, \min\big\{ t \ge 0 \, : \, \B^{(t)}(\A(\S)) \text{ contains a spanning pair for each row of $R$} \big\}.
\]
\end{defn}

\begin{figure}
	\includegraphics[width =  .9 \textwidth]{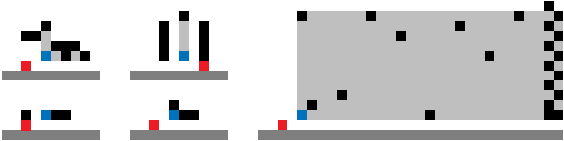}
	\caption{Five configurations (the red and black sites) that do not have a spanning pair for the row above the dark grey row at time $t=0$, but that create a spanning pair (the red and blue sites) at some time $t >0$ by iteration with $\B$. The light grey sites indicate which sites must become infected to create the spanning pair. Note that in each case these sets have minimal cardinality (i.e., if we remove any black site, then iteration of $\B$ will not create the spanning pair).}\label{fig:infector}
\end{figure}

In words, the spanning time $\tau$ is the first time $t$ such that $\B^{(t)}(\A(\S))$ spans all rows of $R$. Since $R$ is up-traversable by $\A(\S)$, it follows by Lemma~\ref{lem:equivuptrav} that $\tau$ must be finite. However, we emphasise that it is possible that $\tau > 0$, see Figure \ref{fig:infector} for some examples.

The central idea in the proof of Lemma~\ref{lem:kvc} is to consider the cases $\tau = 0$ and $\tau > 0$ separately. When $\tau = 0$, the structure is significantly simpler than when $\tau > 0$, which allows for a very sharp estimate. When $\tau > 0$ more complex structures are possible, but more infected sites are required, and this allows us to use a less precise analysis.
\medskip

\subsection{The case $\tau = 0$}
Given a rectangle $R$, let $\F_0(R)$ and $\F_+(R)$ denote the families of all minimal sets $\A \subset R$ such that $R$ is up-traversable by $\A$ and $\tau(R,\A) = 0$ and $\tau(R,\A)>0$, respectively. Let us write $\U_0(R)$ and $\U_+(R)$ for the upsets generated by $\F_0(R)$ and $\F_+(R)$, respectively, i.e., the collections of subsets of $R$ that contain a set $\A \in \F_0(R)$ or $\A \in \F_+(R)$, respectively. 

The following lemma gives a precise estimate of the probability that a rectangle is up-traversable and $\tau=0$.
\begin{lemma}\label{lem:kvctauzero} 
Let $R$ be a rectangle with dimensions $(x,y)$, and let $p \in (0,1)$. Then
\begin{equation}\label{e:nogap}
	\Pp\big( \S \cap R \in \U_0(R) \big) \le (8p^2 x + 8p)^y.	
\end{equation}
\end{lemma}

We will prove Lemma~\ref{lem:kvctauzero} using the first moment method. To be precise, we will show that the expected number of members of $\F_0(R)$ that are contained in $\S$ is at most the right-hand side of~\eqref{e:nogap}. This will follow easily from the following lemma.

\begin{lemma}\label{lem:M0R:size}
Let $R$ be a rectangle with dimensions $(x,y)$, and let $p \in (0,1)$. Then
\[
	|\F_0(R)| \, \le \, \sum_{r=1}^y 8^y {y - 1 \choose r - 1} x^r.
\]
\end{lemma}

\begin{figure}
	\includegraphics[width = .49\textwidth]{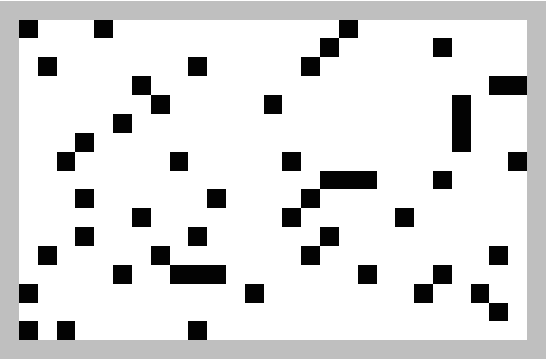}
	\includegraphics[width= .49\textwidth]{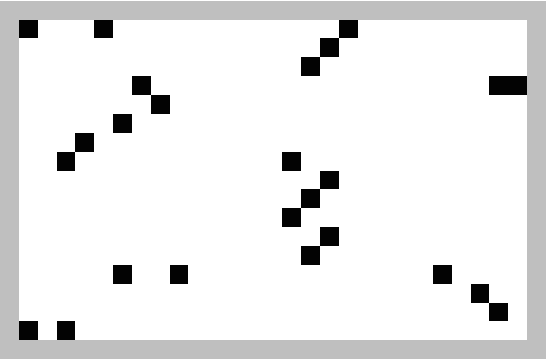}
	\caption{On the left: an up-traversable rectangle. On the right: a minimal set $\A$. Note that $\A$ is sufficient for up-traversability, that $\A$ spans every row (so $\tau=0$), and that $\A$ consists of $8$ paths (so $r=8$).}\label{fig:pairs}
\end{figure}

To count the sets in $\F_0(R)$, we will need to understand their structure. We will show that each set $\A \in \F_0(R)$ can be partitioned into ``paths'' as follows: 

\begin{lemma}\label{lem:A:paths}
Let $R$ be a rectangle with dimensions $(x,y)$, and let $\A \in \F_0(R)$. Then there exists a partition $\A = A_1 \cup \cdots \cup A_r$, where $r = |\A| - y$, with the following property: For each $j \in [r]$, there exists an ordering $(u_1,\ldots,u_{|A_j|})$ of the elements of $A_j$ such that 
\[
	u_i - u_{i-1} \in \big\{ (\pm 2,1), (\pm 1,1) \big\},
\]
for each $2 \le i < |A_j|$, and 
\[
	u_{|A_j|} - u_{|A_j|-1} \in \big\{ (\pm 4,0), (\pm 3,0), (\pm 2,0), (\pm 1,0), (\pm 2,1), (\pm1,1) \big\}.
\]
\end{lemma}
See Figure~\ref{fig:pairs} for an illustration.

\proof
Since $\A$ is a minimal subset of $R$ such that $R$ is up-traversable by $\A$, and $\tau(R,\A) = 0$, it follows from Definition~\ref{def:tau} that $\A$ contains a spanning pair for each row of $R$, and hence (by minimality of $\A$) it follows that $\A$ consists exactly of a union of spanning pairs (one pair for each row) and no other sites. Let these pairs be $\P_1,\ldots,\P_y$, and define a graph on $[y]$ by placing an edge between $i$ and $j$ if $\P_i \cap \P_j$ is non-empty. The sets $A_1,\ldots,A_r$ are simply (the elements of $\A$ corresponding to) the components of this graph.

Let the components of the graph be $C_1,\ldots,C_r$, and note first that each component is a path, since a spanning pair for row $[a,b] \times \{c\}$ is contained in $[a,b] \times [c,c+1]$. Moreover, it follows immediately from this simple fact that if $\P_i \cap \P_j$ is non-empty then $\P_i$ and $\P_j$ must be spanning pairs for adjacent rows (say, $[a,b] \times \{c\}$ and $[a,b] \times \{c+1\}$), and that their common element must lie in $[a,b] \times \{c+1\}$.

Now, consider a component $C_\ell = \{ i_1,\ldots, i_s \}$, set $A_\ell = \bigcup_{j = 1}^s \P_{i_j}$, and note that $|A_\ell| = s + 1$. Let $A_\ell = \{u_1,\ldots,u_{s+1}\}$, and assume (without loss of generality) that $\P_{i_j} = \{u_j,u_{j+1}\}$ for each $j \in \{1,\dots, s\}$. It now follows from the comments above, and the definition of a spanning pair in \eqref{def:spanning:pairs}, that 
\[
	u_i - u_{i-1} \in \big\{ (\pm 2,1), (\pm 1,1) \big\},
\]
for each $2 \le i \le s$, and that 
\[
	u_{s+1} - u_s \in \big\{ (\pm 4,0), (\pm 3,0), (\pm 2,0), (\pm 1,0), (\pm 2,1), (\pm1,1) \big\},
\]
as claimed. Finally, note that $|\A| = y + r$, since $|A_\ell| = |C_\ell| + 1$ for each $\ell \in \{1,\dots, r\}$. \qed
\medskip

\proof[Proof of Lemma~\ref{lem:M0R:size}]
To count the sets $\A \in \F_0(R)$, let us first fix $|\A|$, and the sizes of the sets $A_1,\ldots,A_r$ given by Lemma~\ref{lem:A:paths}. Recall that $r =  |\A| - y$ and that $\A = A_1 \cup \cdots \cup A_r$ is a partition, and note that $|A_j| \ge 2$ for each $j \in \{1,\dots, r\}$, since $A_j$ is a union of spanning pairs. It follows that there are exactly
\[
	{y - 1 \choose r - 1}
\]
ways to choose the sequence $(|A_1|,\ldots,|A_r|)$, where we order the sets $A_j$ so that if $i < j$ then the top row of $A_i$ is no higher than the bottom row of $A_j$. (Note that this is possible because each $A_i$ is a union of spanning pairs for some set of consecutive rows of $R$.) Now, we claim that there are at most
\[
	x \cdot 8^{|A_j| - 1}
\]
ways of choosing the elements of $|A_j|$, given $A_{j-1}$ and $|A_j|$. Indeed, given $A_{j-1}$ we can deduce which is the bottom row of $A_j$, and we have at most $x$ choices for the left-most element $u_1$ of $A_j$ in that row. If $|A_j| = 2$ then (given $u_1$) there are then exactly $8$ choices for the other element $u_2$, since $u_2 - u_1 \in \big\{ (4,0), (3,0), (2,0), (1,0), (\pm 2,1), (\pm1,1) \big\}$. On the other hand, if $|A_j| \ge 3$, then there are at most $4^{|A_j| - 2} \cdot 12 \le 8^{|A_j| - 1}$ choices for the remaining elements of $A_j$ (given $u_1$), by Lemma~\ref{lem:A:paths}, as required.

Now, multiplying together the (conditional) number of choices for each set $A_j$, it follows that 
\[
	|\F_0(R)| \, \le \, \sum_{r=1}^y \sum_{|A_1|,\ldots,|A_r|} \prod_{j = 1}^r \big( x \cdot 8^{|A_j| - 1} \big)  \, \le \, \sum_{r=1}^y 8^y {y - 1 \choose r - 1} x^r,
\]
as claimed, since $\sum_{j = 1}^r (|A_j| - 1) = y$. \qed
\medskip

Lemma~\ref{lem:kvctauzero} now follows by Markov's inequality:

\proof[Proof of Lemma~\ref{lem:kvctauzero}]
Define a random variable $X$ to be the number of sets $\A \in \F_0(R)$ that are entirely infected at time zero, i.e., that are contained in our $p$-random set $\S$. By Markov's inequality and Lemma~\ref{lem:M0R:size}, we have
\begin{align*}
\Pp\big( \S \cap R \in \U_0(R) \big) & \le \Ex_p[X] \le \sum_{r=1}^y 8^y {y - 1 \choose r - 1} x^r p^{y+r} \\
& = \sum_{r=1}^y {y-1 \choose r-1} (8 p^2 x)^r (8 p)^{y-r} = \frac{p x}{1+px} (8p^2 x + 8p)^y
\end{align*}
as required. \qed
\medskip

\subsection{The case $\tau > 0$}

In this section we analyse the event $\S \cap R \in \U_+(R)$. If $R$ is up-traversable by $\S$, then let $\A$ again denote a subset of $\S$ of minimal cardinality such that $R$ is up-traversable by~$\A$. 
By Lemma \ref{lem:equivuptrav} above we know that if $R$ is up-traversable by $\A$, then there must exist a time $t$ at which there is spanning pair in $\B^{(t)}(\A)$ for each row in $R$. The following definition isolates the sites that are responsible for the creation of such spanning pairs.

\begin{defn}\label{def:infector}
Given $\S$ and a row $\ell$, we say that $\M \subset \S$ is an \emph{infector of the row $\ell$} if 
\begin{itemize}
	\item there exists a $t \ge 0$ such that $\B^{(t)}(\M)$ contains a spanning pair for the row $\ell$, and
	\item there does not exist a subset $\M' \subset \M$ such that there exists a $t' \ge 0$ such that $\B^{(t')}(\M')$ contains a spanning pair for the row $\ell$ .
\end{itemize}
We call the bottom-most left-most site in $\M$ the \emph{root} of $\M$. Given $\S$ we write $\Mbb(\S, R)$ for the set of all infectors contained in $\S$ for a row of $R$.  
\end{defn}

Note that spanning pairs are infectors, but that many other configurations are possible: see Figure \ref{fig:infector} for a few examples.

\begin{lemma}[A property of the union of infectors]\label{lem:infector} Suppose $R = [1,x] \times [1,y]$ is up-traversable by $\S$ and that $\A$ is a subset of $\S$ of minimal cardinality with the same property. For each $\ell \in \{1,\dots,y\}$ there exists an infector $\M_\ell$ of row $\ell$ in $\S$ such that
\[
	\bigcup_{\ell=1}^{y} \M_\ell = \A.
\]
\end{lemma}
\proof 
 Let $\A'$ be a subset of $\S$ such that $R$ is up-traversable by $\A'$ and such that $\A'$ is a set with minimal cardinality for this property. By Lemma \ref{lem:equivuptrav}, the event that $R$ is up-traversable by $\A'$ is equivalent to the event that there exists a spanning pair for each row of $R$ after some finite number of iterations of $\A'$ by the bootstrap operator $\B$. This means that for each row $\A'$ contains at least one infector. Note that it is a priori possible that the infectors in $\Mbb(\A', R)$ overlap partially or that an infector for some row $\ell$ is contained in an infector for a row $\ell' \neq \ell$.  Write $(\M^{(i)})_{i=1}^{|\Mbb(\A',R)|}$ for some (arbitrary) ordered list of the infectors, and, for $1 \le s \le |\Mbb(\A', R)|$ write
\[
 	\Mbb^\flat(s) := \bigcup_{i=1}^{s-1} \M^{(i)} \cup \bigcup_{i=s+1}^{|\Mbb(\A',R)|} \M^{(i)}
\]
for the union of the sites of all the infectors except those of $\M^{(s)}$. Now suppose that there exist $1 \le s < t \le |\Mbb(\A', R)|$ such that $\M^{(s)}, \M^{(t)}$ are both infectors of the same row $\ell$ and suppose that $\M^{(s)} \setminus \Mbb^\flat(s) \neq \varnothing$ and $\M^{(t)} \setminus \Mbb^\flat(t) \neq \varnothing$. Then, since $\M^{(s)}$ is an infector for row $\ell$ and the sites in $\M^{(t)} \setminus \Mbb^\flat(t)$ are not needed to create a spanning pair for any other row, $R$ is also up-traversable by the set $\A' \setminus (\M^{(t)} \setminus \Mbb^\flat(t))$, whose cardinality is strictly smaller than $\A'$. This gives a contradiction. Hence, for each row $\ell$ there must exist at most one infector $\M^{(s)}$ with the property that $\M^{(s)} \setminus \Mbb^\flat(s) \neq \varnothing$. Taking their union we obtain $\A$ (i.e., $\A = \A'$).
 \qed
\medskip

Recall that for any set $\Qcal \subset \ZZ^2$ we write $\xb(\Qcal)$ and $\yb(\Qcal)$ for the horizontal and vertical dimensions of that set.
We split the event $\{\S \cap R \in \U_+(R)\}$ according to whether there exists an infector $\M_\ell$ with $\xb(\M_\ell) \ge 6k^2$ or not.

\begin{lemma}[Wide infectors]\label{lem:wideinfector} 
Let $R = [1,x] \times [1,k]$ with $k \ge 3$ such that $ k \ll p^{-1}$, and $x$ such that $ k^5 \ll x \le p^{-2}$, then
\begin{equation}\label{e:wideinfector}
	\Pp\Big(\S \cap R \in \U_+ (R), \max_{\ell=1}^k  \xb(\M_\ell) \ge 6k^2 \Big) = o((8p^2x + 8p)^k).
\end{equation}
\end{lemma}
\proof
Write $\M_{j}$ for the first infector such that $\xb(\M_{j}) \ge 6 k^2$. Since $\M_{j} \subset [1,x] \times [1,k]$, $\yb(\M_{j}) \le k$. Moreover, $\M_{j}$ is the minimal set responsible for the creation of the spanning pair in row $j$, so it must be the case that $\M_{j}$ does not have a gap of more than three consecutive columns. There are at most $xk$ possible positions for the root of the infector. We thus bound \eqref{e:wideinfector} for the range of $x$ and our choice of $k$ from above by
\[
	xk (1-(1-p)^{3k})^{2k^2} \le xk (3pk)^{2k^2} \ll (3 p k)^{6k - 3} 
	\ll \left(800 p^3 k^3 \right)^{k} \ll (8p^2x + 8p)^{k}.\qed
\]

\begin{lemma}[Small infectors]\label{lem:infector2}
	There exist no infectors that are not a single spanning pair that intersect precisely one row, and there exist precisely two infectors that are not a single spanning pair that intersect precisely two rows, up to translations. The cardinality of these infectors is $4$, and they span both rows they intersect.
\end{lemma}
\proof
Let $\M_j$ be the infector for some row $j$.
Write $v$ for an element of the spanning pair for row $j$ that becomes infected due to the bootstrap dynamics on $\M_j$. (It is easy to see that only one element of a spanning pair can arise after time $t=0$, but we do not use this fact.)
Suppose $t$ is the first time such that $\B^{(t)}(\M_j)$ contains a spanning pair. Because $\M_j$ is not a spanning pair, $t \ge 1$. Since $v$ becomes infected at time $t$, it must be the case that $|\Nani(v) \cap \B^{(t-1)} (\M_j)| \ge 3$. Any configuration of three sites in $\Nani(v)$ contains a spanning pair for the row that $v$ is in, so $v$ cannot be in row $j$. By the definition of spanning pairs, \eqref{def:spanning:pairs}, a site can either span the row that it is in, or the row below it, so $v$ is in row $j+1$. 
We conclude that there are no infectors that are not a spanning pair that intersect precisely one row.

By the same argument, if $t \ge 2$, then $\M_j$ must contain a site in row $j+2$, so only infectors that intersect two rows can have $t=1$. 

One can easily verify that the only infectors with $t =1$ that intersect two rows are translations of the configurations $\{(0,0), (0,1), (3,1), (4,1)\}$ and $\{(0,0), (0,1), (-3,1), (-4,1)\}$ (see the configuration in the bottom-left corner of Figure~\ref{fig:infector}). These infectors both have cardinality $4$, and span both rows they intersect.  \qed
\medskip

To analyse $\Pp(\S \cap R \in \U_+(R))$ we again divide $\A$ into the maximal number of disjoint, ``causally independent'' pieces, to which we may apply the BK-inequality. We have seen that when $\tau=0$ these pieces can be described as paths. When $\tau>0$ this is still the case, but now the path structure can be found at the level of the infectors. We partition $\A$ as follows: let $r$ be the largest integer such that there exist sets $B_1, \dots, B_r$ that partition $\A$ (i.e., $B_i \cap B_j = \emptyset$ for all $i \neq j$ and $\A = \bigcup_{i=1}^r B_i$) and such that there exist $r$ pairs of integers $\{(a_i, b_i)\}_{i=1}^r$ such that
\begin{itemize} 
	\item $1 = a_1 \le b_1 \le a_2 \le b_2 \le \dotsm \le a_r \le b_r =k$, and
	\item the event 
		\[
			\{[1,x] \times [a_1, b_1] \text{ is up-trav by }B_1\} \circ \dotsm \circ \{[1,x] \times [a_r, b_r] \text{ is up-trav by } B_r\}
		\]
		occurs.
\end{itemize}

\begin{lemma}[Path structure of $B_1, \dots, B_r$]\label{lem:Bprop}
Let $R = [1,x] \times [1,y]$ and suppose that $R$ is up-traversable by $\S$. Let $\A$ be the subset of $\S$ with minimal cardinality such that $R$ is up-traversable by $\A$. Let $B_1, \dots, B_r$ be the division of $\A$ into disjointly occurring pieces described above. Then the following hold:
\begin{enumerate}
	\item For any row $\ell \in \{1, \dots, y\}$ there exists a unique $i \in \{1,\dots, r\}$ such that $\M_\ell \subseteq B_i$.
	\item If $B_i$ spans rows $\ell, \dots, \ell+m$, then $B_i = \cup_{j=\ell}^{\ell+m} \M_j$.
	\item If $\M_j  \subseteq B_i$ and $j < b_i$, then at least one of the following holds: $\M_j = B_i$;  or there exists a $j' < j$ such that $\M_j \subset \M_{j'}  \subseteq B_i$; or $\M_j \cap \M_{j+1} \neq \emptyset$.
	\item If $\M_j  \subseteq B_i$ and $j = b_i$, then at least one of the following holds: $\M_j = B_i$;  or there exists a $j' < j$ such that $\M_j \subset \M_{j'}  \subseteq B_i$; or $\M_{j-1} \cap \M_{j} \neq \emptyset$.
\end{enumerate}
\end{lemma}

\proof
(a) By construction, $\A = \cup_{i=1}^r B_i$, and $B_i \circ B_j$ occurs if $i \neq j$. By Lemma~\ref{lem:infector}, $\Acal = \cup_{\ell =1}^k \M_\ell$.  Suppose that there exists an $\ell$ such that $\M_\ell\cap B_i\ne \emptyset$ and $\M_\ell\cap B_j\ne \emptyset$ for some $i\ne j$. Without loss of generality, we can further assume that $a_i\le \ell\le b_i$. Since $\M_\ell$ is the minimal set to create a spanning pair for row $\ell$, and that $\M_\ell\cap B_i$ is a strict subset of $\M_\ell$ (since the latter intersects $B_j$, which is disjoint from $B_i$ by assumption), we deduce that $\< \M_\ell\cap B_i\>$ cannot contain a spanning pair for row $\ell$. By Lemma~\ref{lem:equivuptrav}, this means that  $[1,x]\times[a_i,b_i]$ is not up-traversable by $B_i$, which is a contradiction.

(b) By Lemma \ref{lem:infector}, $\A = \cup_{i=1}^k \M_i$. Combined with (a) this gives (b).

(c) Suppose that $B_i$ spans rows $\ell, \dots, \ell+m$ and suppose that there exists a $j < b_i$ such that neither $\M_j = B_i$ nor $\M_j \subset \M_{j'}$ for some $j' < j$, and such that $\M_{j} \cap \M_{j+1} = \emptyset$. Then we can partition 
\[
	B_i = \left(\bigcup_{s=\ell}^j \M_{s} \right) \sqcup \left(\bigcup_{t=j+1}^{\ell+m} \M_{t} \right) =: B_{i,1} \,  \sqcup \, B_{i,2}.
\]
It then follows that
\[
	\{[1,x] \times [\ell, j] \text{ is up-trav by }B_{i,1}\}  \circ \{[1,x] \times [j+1, \ell+m] \text{ is up-trav by } B_{i,2}\}
\]
occurs. This gives a contradiction, since by construction the sets $B_1, \dots, B_r$ are the maximal partition of $\A$ with this property, so such a $j$ does not exist. 
So we conclude that if $\M_j \subset B_i$ but $\M_j \neq B_i$ and $\M_j \nsubseteq \M_{j'}$ for some $j'<j$, then $\M_{j} \cap \M_{j+1} \neq \emptyset$. 

(d) The proof is identical to that of (c), mutatis mutandis.
\qed
\medskip

For all $k,\ell,m,x \in \mathbb{N}$, let $\E_{\ell+1, \ell+m}$ denote the event that a configuration of infected sites $\S$ has the following properties:
\begin{itemize}
	\item $\S \cap ([1,x] \times [\ell +1, \ell + m]) \in \U_+([1,x] \times [\ell +1, \ell + m])$,
	\item the minimal subset $\A$ of $\S$ such that $[1,x] \times [\ell +1, \ell + m]$ is up-traversable by $\A$ cannot be divided into two or more disjointly occurring pieces, i.e., $\A = B_1$ in the construction described above.
	\item $\max_{j=\ell+1}^m \xb(\M_j) < 6k^2$.
\end{itemize}

\begin{lemma}\label{lem:embd} For $k \ge 3$, $\ell + m \le k$ and all $p \in [0,1]$,
\[
	\Pp(\E_{\ell+1, \ell+m}) \le 4p^2 x (12pk^2 + 7p)^{m-1}.
\]
\end{lemma}
\proof
There is at least one infected site in row $\ell +1$, and it can be at $x$ positions.

By Lemma \ref{lem:Bprop}, the event $\E_{\ell+1, \ell+m}$ implies that $\A$ is the union of infectors that are not disjoint.  Since, moreover, none of the infectors are wider than $6k^2-1$, for each of the rows $\ell+2,\dots, \ell+m$ we then need to have at least $1$ infected site in the line-segment $[-6k^2-3,6k^2+3]$ directly above the infected site of the row below it. Finally, row $\ell+m$ must also be spanned, and by Lemma~\ref{lem:infector2} its spanning pair must already be present at time $t=0$, so there must be another infected site in that row, in one of the four positions that can create a spanning pair for line $\ell+m$. We thus bound
\[
	\Pp(\E_{\ell+1, \ell+m})  \le px \cdot (p (12k^2+7))^{m-1} 4p. \qed
\]
\medskip

Write
\[
	\V_{a,b} := \big\{[1,x] \times [a,b] \text{ is up-traversable}\big\}
\]
and
\[
	\V^+_{a,b} := \big\{\S \cap ([1,x] \times [a,b]) \in \U_+([1,x] \times [a,b])\big\} \cap \big\{  \max_{j=a}^b \xb(\M_j) \le 6 k^2 \big\}.
\]
The following lemma states the key inequality for the induction:
\begin{lemma}\label{lem:kvcsplit} For $k \ge 2$,
\[
	\Pp(\V^+_{1,k}) \le \sum_{m=2}^k \sum_{\ell=0}^{k-m} \Pp(\V_{1,\ell}) \Pp(\E_{\ell+1, \ell+m}) \Pp(\V_{\ell + m + 1,k}).
\]
\end{lemma}

\proof Since $\V^+_{1,k}$ occurs, $[1,x] \times [1,k]$ is up-traversable. Let $\A$ be the minimal subset of $\S$ such that $[1,x] \times [1,k]$ is up-traversable with respect to $\A$. Let $B_1, \dots, B_r$ be the subdivision of $\A$ described above. Let $u \in \A$ and $v \in \<\A\> \setminus \A$ be such that $\{u,v\}$ form a spanning pair for the row $i$, while $\A$ does not contain a spanning pair for row $i$. At least one such pair must exist since $\V^+_{1,k}$ occurs. Let $j$ be such that $\M_i \subseteq B_j$ (we can find such a $B_j$ by Lemma~\ref{lem:Bprop}(a)). Suppose that $B_j$ spans exactly the rows $\ell+1, \dots, \ell+m$ (i.e., $a_j = \ell+1$ and $b_j = \ell+m$). Then, by the construction of $B_1, \dots, B_r$ and $\E_{\ell+1, \ell+m}$ we know that 
\[
	\V_{1,\ell} \circ \E_{\ell+1, \ell+m} \circ \V_{\ell+m+1,k}
\]
occurs for $\S$. Applying the BK-inequality and summing over $\ell$ and $m$ gives the asserted inequality.  The sum over $m$ starts at $2$ because by Lemma ~\ref{lem:infector2}, $B_j$ must span at least two rows.\qed
\medskip

\subsection{The proof of Lemma \ref{lem:kvc}}
To begin, assume that $\frac{3 k^2}{p} \le x \le \frac{1}{p^2}$.
We start by proving Lemma \ref{lem:kvc} for the cases where $y \le k$. More precisely, we will prove that
\begin{equation}\label{e:kvcylek}
	\Pr_p(\V_{1,k}) \le \e (8p^2 x +8 p)^k,
\end{equation}
holds for $k \ll p^{-1}$.
We use induction.
The inductive hypothesis is that \eqref{e:kvcylek} holds for $k' \le k-1$ and $ k^5 \ll  x \le p^{-2}$.
To initialise the induction we observe that when $k=1$ there exist four spanning pairs up to translations that intersect one row, so $\Pp(\V_{1,1}) \le 4p^2 x < \e (8p^2x + 8p)$.
When $k=2$ we use Lemma \ref{lem:infector2} to bound
\begin{equation}\label{e:twolines}
	\Pp(\V^+_{1,2}) = 2 x p^4,
\end{equation}
which, combined with Lemma \ref{lem:kvctauzero} yields that	
\[
	\Pp(\V_{1,2}) \le (8p^2x + 8p)^2 + 4xp^4  \le \e (8p^2 x +8 p)^2
\]
when $p$ is sufficiently small.

When $3 \le k \ll p^{-1}$, by \eqref{e:twolines}, Lemmas \ref{lem:kvctauzero}, \ref{lem:wideinfector}, \ref{lem:embd}, and \ref{lem:kvcsplit}, and the induction hypothesis \eqref{e:twolines}, when $p$ is sufficiently small,
\begin{equation}\label{e:V1k}\begin{split}
	\Pr_p(\V_{1,k}) \le & (8p^2x + 8p)^k + \frac{\e-1}{3} (8p^2 x + 8p)^k + 2 x p^4 (k-1) \e^2 (8 p^2 x + 8 p)^{k-2} \\
	& + 4\e^2 x p^2 k \sum_{m=3}^k  (8p^2 x + 8 p)^{k-m+1} (12p k^2 +7p )^{m-1},
\end{split}\end{equation}
where the second term on the right-hand side is due to Lemma~\ref{lem:wideinfector}, and the third and fourth correspond to the $m=2$ and $m\ge 3$ terms in Lemma~\ref{lem:kvcsplit}. 

It is not difficult to show that 
\[
	\sum_{m=3}^k  a^{k-m+1} b^{m-1} = \frac{b^2 a^{k-1} - a b^k}{a-b}.
\]
When $\frac{3 k^2}{p} \le x$ we have $12 p k^2  + 7p \le \tfrac12 (8 p^2 x + 8p)$, so this implies that
\[
	4\e^2 x p^2 k \sum_{m=3}^k  (8p^2 x + 8 p)^{k-m+1} (12p k^2 + 7p )^{m-1} \le  8 \e^2 x p^2 (12 p k^2 + 7p)^2 (8 p^2 x + 8p)^{k-2}.
\]
Inspecting \eqref{e:V1k}, it follows that the desired bound \eqref{e:kvcylek} holds if the following two inequalities hold for $p$ sufficiently small:
\begin{align*}
	2 \e^2 x p^4(k-1)  &< \frac{\e-1}{3}(8 p^2 x + 8 p)^2,\\
	8 \e^2 x p^2 k(12pk^2 + 7p)^2 &< \frac{\e-1}{3} (8p^2 x + 8 p)^2.
\end{align*}
The first inequality holds because $k \ll x$. It is easy to verify that the second inequality holds when $k^5 \ll x \le p^{-2}$.
Substituting the above inequalities into \eqref{e:V1k} proves the claim of Lemma \ref{lem:kvc} for $y \le k$.
\medskip

Now we consider $R = [1,x] \times [1,y]$ for $y$ such that $k < y < x$ (still assuming that $\frac{3 k^2}{p} \le x \le \frac{1}{p^2}$). We cover $R$ with $\lceil y/k \rceil$ rectangles of height~$k$. If $y$ is not divisible by $k$ the covering ``overshoots'': it includes at most $k-1$ rows that are not in $R$. If $R$ is up-traversable, and if the overshoot contains a connected upward path, then all these rectangles are also up-traversable.
The probability that there is a connected path in the overshoot is at least~$p^k$.
It thus follows by the BK-inequality that
\begin{equation}\label{e:bkkvc}
\begin{split}
	\Pp \left(R \text{ is up-trav}\right)  &\le  p^{-k} \prod_{n=1}^{\lceil y /k \rceil} \Pp (\V_{(n-1)k+1, nk})\\
	&
	\le  p^{-k} \e^{y/k} (8p^2 x + 8 p)^{y},
\end{split}
\end{equation}
where these bounds again hold for $p$ sufficiently small. 
This completes the proof of Lemma~\ref{lem:kvc} for the case $\frac{3 k^2}{p} \le x \le \frac{1}{p^2}$. 
\medskip

The case $x < \frac{3 k^2}{p}$ is now easy. Note that if $[1,x] \times [1,y]$ is up-traversable by $\S$, then also $[1,x+a] \times [1,y]$ for any $a \ge 1$ is up-traversable by $\S$ (i.e., up-traversability is a monotone increasing event in the width of the rectangle). Hence, $\Pp([x] \times [y]$ is up-trav$)$ is a monotone increasing function in $x$. The bound thus follows by choosing $x = \frac{3k^2}{p}$ and applying the bound for the case $\frac{3 k^2}{p} \le x \le \frac{1}{p^2}$. 
\qed
\medskip

\section{The probability of simultaneous horizontal and vertical growth}\label{sec:lemWproof}

The lemma below states an upper bound on the probability of an infected rectangle growing both vertically and horizontally, i.e., an upper bound on $\Pp(R \growto R')$ for certain $R \subset R'$.

Let
\[
	\xi := \left\lceil\log^2 \frac1p\right\rceil, \qquad \text{ and } \qquad  \delta_\xi := 1 - 2/\xi.
\]
Recalling Lemma \ref{lem:kvc} and the bound on $f(p,y)$ in Lemma \ref{21} above, let 
\begin{equation}\label{e:gdef}
	\psi(x) :=   -\left(\log(24 p \xi^2 + 8p) + \xi^{-1}\right) \indi_{ \{ x < \frac{3 \xi^2}{ p} \}}-\left(\log(8p^2 x + 8p) + \xi^{-1}\right) \indi_{ \{\frac{3 \xi^2}{ p} \le x \le \frac{1}{p^2}\}},
\end{equation}
and let 
\begin{equation}\label{e:fdef}
	\phi(y) := 
	 \e^{-3py} \indi_{\{y > \frac{4}{p} \log\log \frac1p\}}.
\end{equation}
For two rectangles $R\subset R'$ with dimensions $(x,y)$ and $(x+s,y+t)$, let
\begin{equation}\label{e:Wdef}
	U^p(R,R'):= \delta_\xi \,(t\psi(x+s) + s \phi(y+t)).
\end{equation}
Observe that $\psi$ and $\phi$ are both positive, decreasing, and convex functions (where they are not zero). 

\begin{lemma}\label{lem:W} 
Let $R \subset R'$, with dimensions $(x,y)$ and $(x+s,y+t)$ respectively. Assume that $t \le \frac1p \log^{-4} \frac1p$. Then, for $p$ sufficiently small,
\[
	\Pp (R \growto R') \le  2p^{-\xi} \exp\left(-U^p (R,R') \right).
\]
\end{lemma}

The proof uses a similar strategy as \cite[Proof of Proposition 3.3]{DumEnt13}. Roughly speaking this strategy entails that we ``decorrelate'' the horizontal and vertical growth events needed for $\{R \growto R'\}$.

\proof

If $y+t \le \frac4p \log \log \frac1p$ and $x+s > 1/p^2$, then we use the trivial bound $\Pp(R \growto R') \le 1,$
corresponding to $U^p(R,R')=0$, as required. 

If $y+t \le \frac4p \log \log \frac1p$ and $x+s \le 1/p^2$, then we apply Lemma~\ref{lem:kvc} (with $k = \xi$), again giving the required bound. 

Therefore, we assume henceforth that $y+t > \frac4p \log \log \frac1p$ and $x+s \le \frac{1}{p^2}$.
\medskip

To start, suppose that $(1-\delta_\xi) t \psi(x+s) > \delta_\xi s \phi(y+t)$,
which corresponds to the vertical growth component $t$ being disproportionately large compared to the horizontal growth component $s$.
Then, we can simply ignore the horizontal growth and apply Lemma~\ref{lem:kvc} to bound
\begin{equation}
	\begin{split}
	\Pp(R \growto R') &\le p^{-\xi} \exp(-t \psi(x+s)) \\
	&\le p^{-\xi} \exp(-\delta_\xi t \psi(x+s) - \delta_\xi s \phi(y+t))\\
	& = p^{-\xi} \exp(-U^p(R,R')),
	\end{split}
\end{equation}
and we are done. 
Therefore, let us henceforth also assume that 
\begin{equation}\label{e:tsass}
	(1-\delta_\xi) t \psi(x+s) \le \delta_\xi s \phi(y+t).
\end{equation}
	
\begin{figure}
	\includegraphics[width = 0.85\textwidth]{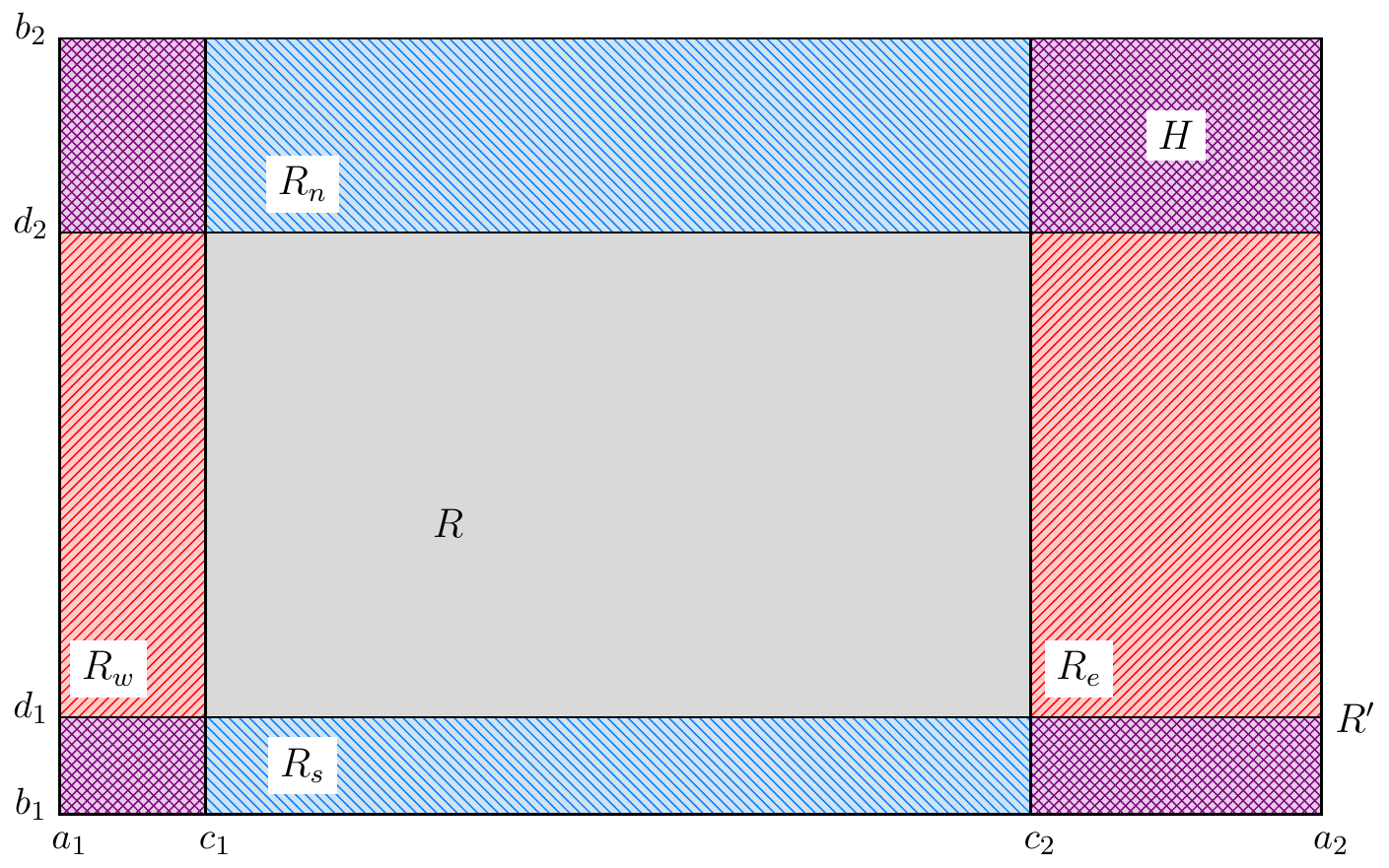}
	\caption{\label{fig:rectangles} The rectangles $R_w$ and $R_e$ are hatched red and the rectangles $R_n$ and $R_s$ are hatched blue. The regions where these four rectangles overlap, collectively called $H$, are cross hatched purple.}
\end{figure}

We identify five (intersecting) regions within the area $R' \setminus R$: the North, South, West, and East regions $R_n$, $R_s$, $R_w$, and $R_e$, and the corner region $H$: for $R'=[a_1,a_2]\times[b_1,b_2]$ and $R=[c_1,c_2]\times[d_1,d_2]$, such that $a_1 \le c_1 < c_2 \le a_2$ and $b_1 \le d_1 < d_2 \le b_2$, we define the sets 
\begin{align*}
	&R_{w}:=[a_1,c_1-1]\times[b_1,b_2]\quad \text{and}\quad R_{e}:=[c_2+1,a_2]\times[b_1,b_2],\\
	&R_{n}:=[a_1,a_2]\times[d_2+1,b_2]\quad \text{and} \quad R_{s}:=[a_1,a_2]\times[b_1,d_1-1],\\
	&H:=R' \setminus \{(x,y):x\in[c_1,c_2]\text{ or }y\in[d_1,d_2]\},
\end{align*}
see Figure~\ref{fig:rectangles}.
Observe that
\[
	\{R \growto R'\} \subset \{R_n \text{ is up-trav}\} \cap \{ R_s \text{ is down-trav}\} \cap \{R_w \text{ and } R_e \text{ are hor-trav}\}.
\]
Let 
\[
	E := \{R_n \text{ is up-trav}\} \cap \{ R_s \text{ is down-trav}\}.
\]
Recall from Definition~\ref{def:infector} above that we write $\Mbb(S,R_n)$ and $\Mbb(\S,R_s)$ for the sets of infectors of $R_n$ and $R_s$ (the latter being a set of infectors suitably defined for down-traversability). By Lemma~\ref{lem:infector}, we are able to determine whether $E$ occurs by inspecting only sites in $\Mbb(S,R_n)$ and $\Mbb(\S,R_s)$. So the event that $R_w$ and $R_e$ are horizontally traversable only depends on $E$ through the information about the intersection of these sets with $H$, the region where the rectangles overlap. Define $\Mbb^\flat_H(\S)$ as the set of all \emph{sites} in $\S \cap H$ contained in either $\Mbb(S,R_n)$ or $\Mbb(\S,R_s)$.
Let $Y$ denote the number of columns in $H$ that contain at least one infected site in $\Mbb^\flat_H(\S)$. We split
\begin{equation}\label{e:Ysplit}
	\Pp(R \growto R') \le \Pp \big(\{R \growto R'\}  \cap \{Y \le s/(2\xi)\} \big) + \Pp(Y > s/(2 \xi)).
\end{equation}

We start by bounding the first term in \eqref{e:Ysplit}. Let $F : = \{Y \le s/(2\xi)\}$. 
We use Lemma~\ref{lem:kvc} with $k = \xi = \lceil \log^2 \frac1p \rceil$ to bound
\begin{equation}\label{e:psifactor}
	\begin{split}
	\Pp( \{R \growto R'\} \cap F) &\le \Pp (R_e \text{ and }R_w \text{ are hor-trav} \mid E \cap F )  \Pp (E )\\
	&\le \Pp (R_e \text{ and }R_w \text{ are hor-trav} \mid E  \cap F ) p^{-\xi} \e^{-t \psi(x+s) }.
\end{split}
\end{equation}

Let $\mfR_n$ denote the set of all sets of $n$ subrectangles of $R_e \cup R_w$ with heights $y+t$, total width $n$, and such that each pair of rectangles in a set $\mfr \in \mfR_n$ are separated by at least one column. I.e., for $\mfr = \{\mfr_i\}_{i=1}^{N(\mfr)} \in \mfR_n$ we have that $\mfr$ is a collection of $N(\mfr)$ strictly disjoint subrectangles with $\sum_{i=1}^{N(\mfr)} \xb(\mfr_i) = n$, and $\yb(\mfr_i) = y+t$ for all $1 \le i \le N(\mfr)$. 
For any $\mfr = \{\mfr_i\}_{i=1}^{N(\mfr)} \in \mfR_n$ define the following two events:
\[ 
	E_1(\mfr) := \left\{\forall 1 \le i \le N(\mfr) : \mfr_i \text{ is horizontally traversable} \right\}
\] 
and
\begin{multline}\label{e:E2def}
	E_2(\mfr) := \bigg\{\Mbb_H^\flat \cap \bigcup_{i=1}^{N(\mfr)} \mfr_i = \varnothing\bigg\} \cap \bigg\{ \nexists \mfr' \in \mfR_{n} : N(\mfr') < N(\mfr) \text{ and }  \Mbb_H^\flat \cap \bigcup_{i=1}^{N(\mfr')} \mfr'_i = \varnothing \bigg\} \\
	\cap \bigg \{\nexists \mfr'' \in \cup_{n'' > n} \mfR_{n''} : \Mbb_H^\flat \cap \bigcup_{i=1}^{N(\mfr'')} \mfr_i = \varnothing \bigg\},
\end{multline}
that is, $E_2(\mfr)$ is the event that $\mfr$ is the partition into the least number of rectangles of total width $n$ that does not intersect $\Mbb_H^\flat$, and that there is no partition of total width greater than $n$ that also does not intersect $\Mbb_H^\flat$.
Observe that
\[	
	\{R_e, R_w \text{ are hor-trav}\} \subseteq  \bigsqcup_{m=0}^s \bigsqcup_{\mfr \in \mfR_{s-m}} (E_1(\mfr) \cap E_2(\mfr)).
\]
Thus,
\begin{multline}\label{e:ReRwsplit}
		\Pp(R_e, R_w \text{ are hor-trav} \mid E \cap F) \\
		 \le \sum_{m=0}^{s /(2\xi)} \sum_{\mfr \in \mfR_{s-m}} \Pp(E_1(\mfr) \mid E_2(\mfr) \cap E \cap F) \Pp(E_2(\mfr) \mid E \cap F),
\end{multline}
where we used that the sum may be restricted to $m \le s/(2\xi)$ by the conditioning on $F$.		
Now note that the events $E$ and $F$ can be verified by inspecting only $\Mbb^\flat_H$, which, on the event $E_2(\mfr)$ is contained in $H \setminus \mfr$, while $E_1(\mfr)$ by definition only depends on the sites in $\mfr$, so that conditionally on $E_2(\mfr)$ the event $E_1(\mfr)$ is independent of $E$ and $F$. We may thus write
\[\begin{split}
	  \Pp(E_1(\mfr) \mid E_2(\mfr) \cap E \cap F) &= \frac{\Pp(E_1(\mfr) \cap E \cap F \mid E_2(\mfr))}{\Pp(E \cap F \mid E_2(\mfr))} \\
	  &= \frac{\Pp(E_1(\mfr) \mid E_2(\mfr)) \Pp(E \cap F \mid E_2(\mfr))}{\Pp(E \cap F \mid E_2(\mfr))} = \Pp(E_1(\mfr) \mid E_2(\mfr)).
\end{split}\]

Observe that for any fixed $\mfr$ the event $E_1(\mfr)$ is increasing. Indeed, adding more sites to $\S$ can either make horizontal traversal occur when it did not before, or else, have no effect. We claim that the event $E_2(\mfr)$, on the other hand, is the intersection of three decreasing events, and hence itself a decreasing event. To see this, observe that the first event in \eqref{e:E2def} is decreasing because adding more sites to $\S$ cannot decrease the total width of $\Mbb_H^\flat$, since it is the union of \emph{all} infectors intersecting $H$ (not only those of minimal cardinality for a given row). The second event in \eqref{e:E2def} is likewise decreasing, because increasing $\Mbb_H^\flat$ cannot decrease the minimal number of rectangles of a partition that does not intersect $\Mbb_H^\flat$, unless it also decreases the total width of that partition. The third event is decreasing because increasing $\Mbb_H^\flat$ cannot decrease its total width.
Therefore, we may apply the FKG-inequality to obtain
\[
	\Pp(E_1(\mfr) \mid E_2(\mfr)) \le \Pp(E_1(\mfr)),
\]
and we may thus further bound the right-hand side of \eqref{e:ReRwsplit} by
\[
	\sum_{m=1}^{s /(2\xi)} \sum_{\mfr \in \mfR_{s-m}} \Pp(E_1(\mfr)) \Pp(E_2(\mfr) \mid E \cap F).
\]
Uniformly for any fixed $\mfr \in \mfR_{s-m}$ with $m \le s/(2\xi)$, by Lemma~\ref{21},
\[\begin{split}
	\Pp(E_1(\mfr)) &= \prod_{i=1}^{N(\mfr)} \Pp(\mfr_i \text{ is hor-trav}) \le \exp\left(-\sum_{i=1}^{N(\mfr)} (\xb(\mfr_i) -2)f(p,y+t)\right)\\
	& \le \exp\left(-s(1-1/\xi) f(p,y+t) \right)\\
	& \le \exp\left(- \delta_\xi s \phi(y+t) \right),
	\end{split}
\]
where the final inequality follows from Lemma~\ref{21}(b) when $p$ is sufficiently small.
Inserting this bound in \eqref{e:ReRwsplit}, we proceed by using that the events $E_2(\mfr)$ are mutually disjoint for all $\mfr$ to bound
\begin{equation}\label{e:narrowomega}
	\begin{split}
	\Pp(R_e, R_w \text{ are hor-trav} \mid E \cap F)  &\le  \exp\left(- \delta_\xi s \phi(y+t)\right) \sum_{m=1}^{s /(2\xi)} \sum_{\mfr \in \mfR_{s-m}} \Pp(E_2(\mfr) \mid E \cap F)\\
	& \le \exp\left(- \delta_\xi s \phi(y+t)\right).
\end{split}
\end{equation}
Combining \eqref{e:psifactor} and \eqref{e:narrowomega} we bound the first term in \eqref{e:Ysplit} by $ p^{-\xi} \exp(-U^p(R,R'))$.
\medskip

Now we bound the second term in \eqref{e:Ysplit}. 
If $Y > s/(2\xi)$ then at least $s/(2 \xi)$ out of $s$ columns are non-empty. The probability that a column is non-empty is $1-(1-p)^t \le 2pt$ (when $p$ is sufficiently small). Therefore, $\Pr(Y > s/(2\xi)) \le \Pr($Bin$(s, 2pt) > s/(2\xi))$. We use Chernoff's bound that $\Pr($Bin$(n, p) > q) \le \e^{-q}$ when $q>np$ 
to estimate 
\[
	\Pp(Y > s/(2\xi)) \le \exp(-s/(2\xi))
\]
(here we used that $t \le \frac1p \log^{-4} \frac1p$).
Observe that since $\xi = \lceil \log^2 \frac1p \rceil$, $\delta_\xi = 1-\xi^{-1}$, and $\phi(y+t) > \log^{-12} \frac1p$ by our assumption that $y+t > \frac4p \log \log \frac1p$, we have
\[
	\exp(-s/(2 \xi)) \le \exp\left(-\frac{\delta_\xi}{1-\delta_\xi} s \phi(y+t)\right).
\]
Now recall our assumption \eqref{e:tsass} that $(1-\delta_\xi) t \psi(x+s) \le \delta_\xi s \phi(y+t)$.
Applying this inequality twice, it follows that
\[
	 \frac{\delta_\xi}{1-\delta_\xi} s \phi(y+t) \ge t \psi(x+s) \ge \delta_x t \psi(x+s) + \delta_\xi s \phi(y+t).
\]
We thus have $\Pp(Y > s/(2\xi)) \le \exp(-U^p(R,R'))$, as required.

Applying the bounds for the two cases to \eqref{e:Ysplit} completes the proof (using the crude upper bound $p^{-\xi} + 1 \le 2 p^{- \xi}$ for $p$ sufficiently small).
\qed

\section{The upper bound of Theorem \ref{thm:IS}}\label{sec:upper}
\begin{prop}\label{lower}
Let $p > 0$ and $\frac{1}{3p} \log \frac{1}{p} \le y \le \frac{1}{p} \log \frac{1}{p}$ and $\frac{1}{p^2} \le x \le \frac{1}{p^5}$. Then \emph{
\[
	\Pp([x] \times [y] \text{ is IF})  \le  \exp\left( - \frac{2C_1}{p} \log^2 \frac{1}{p}  + (2 C_2 + o(1)) \frac{1}{p} \log \frac{1}{p} \right).
\]}
\end{prop}

\subsection{Notation and definitions}
Before we proceed with the proof, we must introduce some more notation and a few definitions. 
Our proof uses \emph{hierarchies}. The notion of hierarchies is due to Holroyd \cite{Hol03}, and is common to much of the bootstrap percolation literature since. Here we use a definition of a hierarchy that is similar to the one in \cite{DumEnt13}:
\begin{defn}[Hierarchies]\color{white}.\color{black}
\begin{enumerate}
	\item \textsc{Hierarchy, seed, normal vertex, and splitter:} A \emph{hierarchy} $\HH$ is a rooted tree with out-degrees at most three\footnote{In the original construction of a hierarchy by Holroyd \cite{Hol03} for the standard model, hierarchies have out-degree at most two. The fact that we need out-degree three corresponds to the fact that the anisotropic model requires three infected sites in a neighbourhood. As a result, it is possible that a rectangle is internally filled by a set of three but not two disjoint internally filled smaller rectangles. Our definition of hierarchies reflects this. See also \cite{DumEnt13}.} and with each vertex $v$ labeled by non-empty rectangle $R_v$ such that $R_v$ contains all the rectangles that label the descendants of~$v$. If the number of descendants of a vertex is $0$, we call the vertex a \emph{seed.}\footnote{Note that although similar, this definition of a seed is different than the one used in the previous section.} If the vertex has one descendant, we call it a \emph{normal} vertex, and we write $u\mapsto v$ to indicate that $u$ is a normal vertex with (unique) descendant $v$. If the vertex has two or more descendants, we call it a \emph{splitter} vertex. We write $N(\HH)$ for the number of vertices in the tree $\HH$.
	\item \textsc{Precision:} A hierarchy \emph{of precision} $Z$ (with $Z \ge  1$) is a hierarchy that satisfies the following conditions:
	\begin{itemize}
		\item[(1)] If $w$ is a seed, then $\xb(R_w) \ge 2$ and $\yb(R_w)<2Z$, while if $u$ is a normal vertex or a splitter, then $\yb(R_u)\ge 2Z$.
		\item[(2)] If $u$ is a normal vertex with descendant $v$, then $\yb(R_u)-\yb(R_v) \le 2Z$.
		\item[(3)] If $u$ is a normal vertex with descendant $v$ and $v$ is either a seed or a normal vertex, then $\yb(R_u)-\yb(R_v) >Z$.
		\item[(4)] If $u$ is a splitter with descendants $v_1,\dots,v_i$ and $i \in \{2,3\}$, then there exists $j\in\{1,\dots,i\}$ such that $\yb(R_{u})-\yb(R_{v_j})> Z.$ 
	\end{itemize}
	\item\textsc{Presence:} Given a set of infected sites $\S$ we say that a hierarchy $\HH$ is \emph{present in $\S$} if all of the following events occur \emph{disjointly}:
	\begin{itemize}
		\item[(1)] For each seed $w$, $R_w = \langle R_w \cap \S\rangle$ (i.e., $R_w$ is internally filled by $\S$).
		\item[(2)] For each normal $u$ and every $v$ such that $u \mapsto v$, $R_u = \langle (R_v \cup \S) \cap R_u \rangle$ (i.e., the event $\{R_v  \growto R_u\}$ occurs on $\S$).
		
	\end{itemize}
	\item\textsc{Goodness:} Similar to \cite{GraHolMor12}, we say that a seed $w$ is \emph{large} if $Z/3 \le \yb(R_w) \le Z$. We call a hierarchy \emph{good} if it has at most $\log^{11} \frac1p$ large seeds, and we call it \emph{bad} otherwise.
\end{enumerate}
\end{defn}

\subsection{Outline of the proof of Proposition \ref{lower}}
In this section we give the proof of Proposition \ref{lower} subject to Lemma  \ref{lem:Wbd} below. We prove Lemma \ref{lem:Wbd} in Section \ref{sec:Wbd}.
\medskip

Let $\HH_{Z,R}$ denote a hierarchy with root $R$ and precision $Z$. Let $\Hb_{Z,R}$ denote the set of all $\HH_{Z,R}$. Likewise, let $\mathbb{H}_{Z,R}^{\sss \mathrm{good}}$ and $\mathbb{H}_{Z,R}^{\sss \mathrm{bad}}$ denote the subsets of good and bad hierarchies in $\Hb_{Z,R}$. Lastly, given a set of hierarchies $\mathbb{H}$ and a rectangle $R$, define the event
\[
	\mathcal{X}(R; \mathbb{H}) := \{\S \in \{0,1\}^{\ZZ^2} \, : \, \exists \HH \in \mathbb{H} \text{ such that } \HH \text{ is present in }\S \cap R\}.
\]
\begin{lemma}\label{lem:exist} Let $R$ be a rectangle with $\xb(R) \ge 2$ and let $Z \ge 3$. If $R$ is internally filled, then there exists a hierarchy $\HH_{Z,R} \in \Hb_{Z,R}$ that is present, i.e., $\mathcal{X}(R; \mathbb{H}_{Z,R})$ occurs.
\end{lemma}

The proof of this lemma is the same as the proof of \cite[Proposition 3.8]{DumEnt13} so we do not repeat it here. (But note that it does not matter that our definition of hierarchies uses ``internally filled'' rather than ``$k$-occurs''.)\medskip

Throughout this paper, let
\[
	Z_p :=\frac 1p \log^{-8} \frac{1}{p}. 
\]
Conform the hypothesis of Proposition \ref{lower} we restrict ourselves to hierarchies with root label $R_p$ of dimensions $(x,y)$ such that
\[
	\frac{1}{p^2} \le x  \le \frac{1}{p^5} \qquad \text{ and } \qquad \frac{1}{3p} \log \frac1p \le y \le \frac{1}{p} \log \frac1p.
\]
For the sake of simplicity we often suppress subscripts $Z_p$ and $R_p$.

We bound the good and bad hierarchies separately:
\begin{equation}\label{e:Rkcrossbd}
	\Pp(R_p \text{ is IF}) \le 
	\Pp(\mathcal{X}(R_p; \mathbb{H}^{\sss \mathrm{good}}))
	+ \Pp(\mathcal{X}(R_p; \mathbb{H}^{\sss \mathrm{bad}})).
\end{equation}

We bound the second term with the following lemma:
\begin{lemma}\label{lem:manyseeds} As $p$ tends to $0$ we have \emph{
\[
	 \Pp(\mathcal{X}(R_p; \mathbb{H}^{\sss \mathrm{bad}}) ) \le \exp\left(- \Omega\left(\frac1p \log^3 \frac1p \right)\right).
\]}
\end{lemma}
\proof
We claim that if $R$ is a large seed, i.e., $Z_p/3 \le \yb(R) \le Z_p$, then
\[
	\Pp(R \text{ is IF}) \le  \exp\left(-\Omega\left(\frac1p \log^{-8} \frac1p\right)\right).
\]
To see that this is indeed the case we consider the cases $x \ge 1/p$ and $x < 1/p$ separately. 
For the case $x \ge 1/p$, the bound follows from Lemma \ref{21}(c):
\[
	\begin{split}
	\Pp(R \text{ is IF}) & \le \Pp(R \text{ is hor-trav}) \le \exp\left(-(x-2)\left(\tfrac12 p y - 3 p^2 y^2 \right)\right)\\
		& \le \exp\left(-\left(\frac1p-2\right)\left(\tfrac12 p \cdot \frac1{3p} \log^{-8} \frac1p - 3 p^2 \cdot \frac{1}{p^2} \log^{-16} \frac1p\right)\right)\\
		& \le \exp\left(-\Omega\left(\frac{1}{p} \log^{-8} \frac1p \right)\right).
	\end{split}
\]
For the case $x < 1/p$, the bound follows from Lemma~\ref{lem:kvc} with $k=2$ and $p$ sufficiently small:
\[
	\begin{split}
		\Pp(R \text{ is IF}) & \le \Pp(R \text{ is up-trav}) \le p^{-2} \e^{y/2} (104 p)^y\\
		& = \exp\left(2 \log \frac1p + \frac{y}2 + y \log(108) - y \log\frac1p \right)\\
		& = \exp\left(- (1+o(1)) y \log \frac1p	\right) = \exp \left(-\Omega\left(\frac{1}{p} \log^{-8} \frac1p \right)\right).
	\end{split}
\]

Now consider the event $\mathcal{X}(R_p; \mathbb{H}^{\sss \mathrm{bad}})$. This event implies that there exists a hierarchy $\HH$ that is present and bad, which by definition means that more than $\log^{11} \frac{1}{p}$ rectangles of size between $Z_p/3$ and $Z_p$ are internally filled disjointly. Since $R_p$ contains at most $\frac{1}{p^6} \log \frac1p$ sites, the probability of this event is smaller than
\[
	\left(\frac{1}{p^6} \log \frac1p \cdot \e^{-\frac{c}{p} \log^{-8} \frac1p }\right)^{ \log^{11} \frac{1}{p}} \le \exp\left(- \Omega\left(\frac1p \log^{3} \frac1p\right)\right).\qed
\]
\smallskip

We bound the first term of \eqref{e:Rkcrossbd} as follows:
\begin{equation}\label{e:uniformX}
	\Pp(\mathcal{X}(R_p; \mathbb{H}^{\sss \mathrm{good}}))  \le |\mathbb{H}^{\sss \mathrm{good}}|  \max_{\HH \in \mathbb{H}^{\sss \mathrm{good}}} \Pp(\HH \text{ is present}).
\end{equation}

Now we apply the following lemma.
\begin{lemma}\label{lem:hierarchy}
The number of good hierarchies satisfies
\[
	|\mathbb{H}^{\sss \mathrm{good}}| \le \e^{O(1/p)}.
\]
\end{lemma}
\proof
We start by observing that any good hierarchy $\HH \in \mathbb{H}^{\sss \mathrm{good}}$ has root $R_p$ such that $x \le \frac1{p^5}$ and $\yb(R_p) \le \frac1p \log \frac1p$ and precision $Z_p = \frac1p \log^{-8} \frac1p$, so its height $h(\HH)$ is bounded from above by 
\[
	h(\HH) \le \frac{\yb(R_p)}{Z_p}   \le \log^9 \frac{1}{p}.
\]
Moreover, since there are at most $\log^{11} \frac{1}{p}$ large seeds in a good hierarchy, the number of vertices $N_\HH$ in the hierarchy $\HH$ obeys
\begin{equation}\label{e:vertices}
	N_\HH \le \log^{11} \frac{1}{p} \cdot \log^9 \frac{1}{p} = \log^{20} \frac{1}{p}.
\end{equation}
Each vertex of a hierarchy has 0, 1, 2 or 3 descendants, so there are at most $4^{\log^{20} \frac{1}{p}}$ unlabelled trees corresponding to the good hierarchies. Finally, since each vertex of a hierarchy is labelled by a sub-rectangle of $R_p$ with $\xb(R_p) \le p^{-5}$, the number of choices for each label is bounded from above by 
\[
	\xb(R_p)^2 \cdot \yb(R_p)^2 \le p^{-10} \cdot \frac{1}{p^2} \log^2 \frac1p \le p^{-13},
\]
so
\[
	|\mathbb{H}^{\sss \mathrm{good}}| \le \left(4p^{-13}\right)^{\log^{20} \frac{1}{p}} = \e^{13 \log^{21} \frac1p + \log 4 \log^{20} \frac1p} \le \e^{O(1/p)}. \qed
\]
\smallskip

By Lemma \ref{lem:hierarchy} it suffices to give a uniform bound on the probability that a given hierarchy is present, if the hierarchy is good. 
Indeed, it remains to show that 
\[
	\max_{\HH \in \mathbb{H}^{\sss \mathrm{good}}} \Pp(\HH \text{ is present}) \le \exp\left( -  \frac{2C_1}{p} \log^2 \frac{1}{p}  + (2 C_2 + o(1)) \frac{1}{p} \log \frac{1}{p} \right).
\]

Before we proceed, let us deal with a small technical issue: the possibility of ``wide'' seeds (Lemma \ref{lem:Wbd} below does not work in their presence).
Observe that if the hierarchy $\HH$ that maximises the probability of being present contains a seed with label $R_s$ such that $\xb(R_s) > \frac{1}{3p} \log^{12} \frac{1}{p}$ (i.e., the seed is extremely wide), then the probability that $\HH$ is present is bounded by the probability that $R_s$ is horizontally traversable, which, by Lemma \ref{21}(c) can be bounded as follows:
\[\begin{split}
	\Pp(R_s \text{ is IF}) & \le \Pp(R_s \text{ is hor-trav})\\
	& \le \exp\left(-(\xb(R_s)-2) \left(\tfrac12 p \yb(R_s) - 3 p^2 \yb(R_s)^2\right)\right)\\
	& \le \exp\left(-\left(\frac1p \log^{12} \frac1p -2\right)\left(\tfrac12 \log^{-8} \frac1p - 3 \log^{-16} \frac1p \right)\right)\\
	& = \exp\left(-\Omega\left(\frac1p \log^3 \frac1p \right)\right),
	\end{split}
\]
where for the third inequality we used the assumption on $\xb(R_s)$ and that $\yb(R_s) < 2Z = \frac2p \log^{-8} \frac 1p$.
Proposition \ref{lower} thus holds for hierarchies with wide seeds.
Let us therefore assume from here on that $\xb(R_s) \le \frac{1}{3p} \log^{12} \frac{1}{p}$ for all seeds.

By the BK-inequality we have
\begin{equation}\label{e:prodH}
	\Pp(\HH \text{ is present}) \le \prod_{u \text{ seed}} \Pp (R_u \text{ is IF}) \prod_{v \mapsto w } \Pp \left(R_v \growto R_w \right)
\end{equation}
(we ignore here the contributions from splitter vertices).

The following lemma is used to determine a bound for the product of the seeds:
\begin{lemma}\label{lem:simpler lemma}
Given a hierarchy $\HH_{Z,R}$, let $N_{\rm seed}$ denote the number of seeds of the hierarchy $\HH_{Z,R}$, and let $u_1,\dots,u_{N_{\rm seed}}$ be an arbitrary ordering of the seeds of $\HH_{Z,R}$. Then \emph{
\[
\prod_{u\text{ seed}}\Pp(R_u \text{ is IF}) \leq \prod_{n=1}^{N_{\rm seed}}\Pp \left(\tilde{R}_n \growto \tilde{R}_{n+1}\right)
\]}
where
\[
	\tilde{R}_n= \left[1,\xb(R_{u_1})+\dots+\xb(R_{u_n}) \right]\times \left[1,\yb(R_{u_1})+\dots+\yb(R_{u_n}) \right].
\]
\end{lemma}
\proof For any $p>0$ and any rectangle $R'$ with dimensions $(x,y)$ with $\min\{x,y\} \ge 2$ and any $a \ge 2$, $b \ge 1$, we have
\begin{equation}\label{e:simpleeq}
	\Pp(R' \text{ is IF}) \le \Pp([1,a] \times [1,b] \growto [1,a+x] \times [1,b+y]).
\end{equation}
Indeed, if the rectangle $[a+1,a+x] \times [b+1, b+y]$ is internally filled, then the event $\{[1,a] \times [1,b] \growto [1,a+x] \times [1, b+y]\}$ occurs. An application of the FKG-inequality thus gives~\eqref{e:simpleeq}. 

Any seed of a hierarchy must have dimensions at least $(2,1)$ by definition, so an iterated application of \eqref{e:simpleeq} completes the proof.\qed
\medskip

Recall the definition of $U^p(R,R')$ in \eqref{e:Wdef} above. We use Lemmas  \ref{lem:W} and \ref{lem:simpler lemma} to bound the first product on the right-hand side of \eqref{e:prodH}:
\[\begin{split}
	\prod_{u \text{ seed}} \Pp (R_u \text{ is IF}) &\le \prod_{n=1}^{N_{\rm seed}}\Pp(\tilde{R}_n \growto \tilde{R}_{n+1})\\
	& \le 2p^{-\xi N_{\rm seed} } \exp\left(-\sum_{n=0}^{N_{\rm seed}} U^p (\tilde R_{n} , \tilde R_{n+1})\right),
\end{split}\]
where $\tilde R_0 = [1] \times [1]$.

To bound the second product of \eqref{e:prodH}, we use the following lemma:
\begin{lemma}\label{less simple lemma}
Let $p>0$. Let $N_{\rm splitter}$ denote the number of splitter vertices of the hierarchy~$\HH$. Then there exists an integer $\hat N = \hat N(\HH) \ge 1$ and a sequence of nested rectangles $\hat{R}_0\subset \dotsm \subset \hat R_{\hat N}$ with the following properties:
\begin{itemize}
\item $\hat R_0 = \tilde R_{N_{\rm seed}}$ (with $\tilde R_{N_{\rm seed}}$ as defined in Lemma \ref{lem:simpler lemma} above),
\item $\hat R_{\hat N}$ has dimensions larger than $R$,
\item $\yb(\hat R_{n+1})-\yb(\hat R_n)\le \frac{1}{p} \log^{-8} \frac1p$ for every $0\le n\le \hat N-1$,
\item for $p$ sufficiently small,
\[
	\prod_{v\mapsto w}\Pp(R_w \growto R_v)\le 2p^{-\xi N_{\rm splitter}}\prod_{n=0}^{\hat N-1}\exp \left(-U^p(\hat R_n, \hat R_{n+1})\right).
\]
\end{itemize} 
\end{lemma}
The proof of this lemma goes by induction, using Lemma \ref{lem:W}, and it is essentially the same as the proof of \cite[Lemma 3.11]{DumEnt13}, so we omit it here.
\medskip

We use Lemma \ref{less simple lemma} to determine that there exist rectangles $\hat R_1\subset\cdots\subset \hat R_{\hat N}$ satisfying the conditions of the lemma such that
\[
	\prod_{v\mapsto w}\Pp(R_v \growto R_w) \le 2p^{-\xi N_{\rm splitter}}\prod_{n=0}^{\hat N-1}\exp \left(-U^p(\hat R_n,\hat R_{n+1})\right).
\]
Using Lemmas \ref{lem:W} and \ref{less simple lemma}
and writing $(R_n)_{n=0}^{N}$ for the concatenation of the sequences $(\tilde R_n)_{n=0}^{N_{\rm seed}}$ and $(\hat R_n)_{n=0}^{\hat N}$, i.e.,
\[
	(R_n)_{n=0}^N := (\tilde R_0, \tilde R_1, \dots, \tilde R_{N_{\rm seed}}, \hat R_1, \dots, \hat R_{\hat N})
\]
 with $N := N_{\rm seed} + \hat N$,
 we bound
\begin{equation}\label{e:Hoccbd}
	\Pp (\HH \text{ is present}) \le   4p^{-\xi (N_{\rm seed} + N_{\rm splitter})}\, \exp\left(-\sum_{n=0}^{N -1} U^p(R_n, R_{n+1}) \right).
\end{equation}

To bound the first factor in \eqref{e:Hoccbd} we use the following lemma:
\begin{lemma}\label{lem:seedsplitter}
Any good hierarchy satisfies
\[
	4p^{-\xi (N_{\rm seed} + N_{\rm splitter})} \le \e^{O(1/p)}.
\]
\end{lemma}
\proof
By \eqref{e:vertices} there are at most $\log^{20} \frac1p$ vertices in a good hierarchy, and $\xi = \lceil\log^2 \frac1p \rceil$, so for any $\HH \in \mathbb{H}^{\sss \mathrm{good}}$,
\[
	4p^{-\xi (N_{\rm seed} + N_{\rm splitter})} \le 4p^{-\lceil \log \frac1p \rceil^{22}} = 4\left(\e^{\log \frac1p}\right)^{ \lceil \log \frac1p \rceil^{22}} \le \e^{O(1/p)}. \qed
\]
\smallskip

The final ingredient of the proof is the following lemma:
\begin{lemma}\label{lem:Wbd}
Let $\mathcal{R}_N = \{R_n\}_{n=0}^N$ be a sequence of increasing, nested rectangles such that
\begin{itemize}
	\item $R_{0} = [1]\times[1]$, 
	\item $y_1 \le \frac{2}{p} \log^{-8} \frac1p$ and $\frac{1}{3p} \log \frac1p \le y_N \le \frac1p \log \frac1p$,
	\item $x_1 \le \frac{1}{3p} \log^{12} \frac{1}{p}$ and $x_N \ge \frac{1}{p^{2}}$.
\end{itemize}
Then
\[
	\exp\left( -\sum_{n=0}^{N-1} U^p(R_n, R_{n+1})\right) \le \exp\left(- \frac{2C_1}{p} \log^2 \frac{1}{p}  + (2 C_2 + o(1)) \frac{1}{p} \log \frac{1}{p} \right).
\]
\end{lemma}
The proof involves a longer computation, so we defer it to Section \ref{sec:Wbd}.
\medskip

\proof[Proof of Proposition \ref{lower} subject to Lemma \ref{lem:Wbd}]
We combine the above lemmas and the bounds derived in the discussion to conclude that
\[\begin{split}
	\Pp(R \text{ is IF}) &\le \Pp(\mathcal{X}(R_p; \mathbb{H}^{\sss \mathrm{good}}))
	+ \Pp(\mathcal{X}(R_p; \mathbb{H}^{\sss \mathrm{bad}})) \qquad \text{{\tiny [Lemma \ref{lem:exist} \& \eqref{e:Rkcrossbd}]}} \\
	\text{{\tiny [\eqref{e:uniformX} \& Lemma \ref{lem:manyseeds}]}} \quad & \le |\mathbb{H}^{\sss \mathrm{good}}|  \max_{\HH \in \mathbb{H}^{\sss \mathrm{good}}} \Pp(\HH \text{ is present}) + \e^{-\Omega\left(\frac1p \log^3 \frac1p \right)}\\
	\text{{\tiny [\eqref{e:prodH} \& Lemma \ref{lem:hierarchy}]}}\quad &\le \e^{O(1/p)} \prod_{u \text{ seed}} \Pp (R_u \text{ is IF}) \prod_{v \mapsto w } \Pp \left(R_v \growto R_w \right) + \e^{-\Omega\left(\frac1p \log^3 \frac1p \right)}\\
	\text{{\tiny [Lemma \ref{lem:simpler lemma}]}}\quad &\le \e^{O(1/p)} \prod_{n=1}^{N_{\rm seed}}\Pp \left(\tilde{R}_n \growto \tilde{R}_{n+1}\right)\prod_{v \mapsto w } \Pp \left(R_v \growto R_w \right) + \e^{-\Omega\left(\frac1p \log^3 \frac1p \right)}\\
	\text{{\tiny [Lemmas \ref{lem:W} \& \ref{less simple lemma}]}}\quad &\le \e^{O(1/p)} 2p^{-\xi N_{\rm seed} } \exp\left(-\sum_{n=0}^{N_{\rm seed}} U^p (\tilde R_{n} , \tilde R_{n+1})\right) \\
	& \quad \times 2p^{-\xi N_{\rm splitter}}\exp \left(-\sum_{n=0}^{\hat N-1} U^p(\hat R_n,\hat R_{n+1})\right) + \e^{-\Omega\left(\frac1p \log^3 \frac1p \right)}\\
	\text{{\tiny [Lemma \ref{lem:seedsplitter}]}}\quad&\le \e^{O(1/p)} \exp\left( -\sum_{n=0}^{N-1} U^p(R_n, R_{n+1})\right)+ \e^{-\Omega\left(\frac1p \log^3 \frac1p \right)}\\
\text{{\tiny [Lemma \ref{lem:Wbd}]}}\quad	& \le \exp\left(- \frac{2C_1}{p} \log^2 \frac{1}{p}  + (2 C_2 + o(1)) \frac{1}{p} \log \frac{1}{p} \right),
\end{split}
\]
as claimed. \qed
\medskip

It remains to prove Lemma \ref{lem:Wbd}. We do this in the upcoming section.


\section{Variational principles: proof of Lemma \ref{lem:Wbd}}\label{sec:Wbd}
To prove Lemma \ref{lem:Wbd} we will start by setting up some variational principles, similar to \cite[Section 6]{Hol03}. We start with a few general lemmas.

Assume throughout this section that $f(x)$ and $g(y)$ are positive, non-increasing, convex, Riemann-integrable functions. 
Let $\RR_+ = (0,\infty)$ and for $\ul a=(a_1,a_2) \in \RR_+^2$ and $\ul b = (b_1, b_2) \in \RR_+^2$, write $\ul a \le \ul b$ if $a_1 \le b_1$ and $a_2 \le b_2$. For $\ul a, \ul b \in \mathbb{R}^2_{+}$ with $\ul a \le \ul b$ and any path $\gamma$ from $\ul a$ to $\ul b$, define
\begin{equation}
	w_{f,g}(\gamma) := \int\limits_\gamma f(x) \d y + g(y) \d x,
\end{equation}
and
\begin{equation}\label{e:Windef}
	W_{f,g}(\ul a, \ul b) := \inf_{\gamma \,:\, \ul a \to \ul b} \int\limits_\gamma f(x) \d y + g(y) \d x.
\end{equation}

To start, an elementary lemma:
\begin{lemma}\label{lem:12} If $\ul a \le \ul b \le \ul c$, then $W_{f,g} (\ul a, \ul b) + W_{f,g}(\ul b, \ul c) \ge W_{f,g} (\ul a, \ul c)$.
\end{lemma}
The proof is easy (see \cite[Section~6]{Hol03}).
\medskip

Let
\begin{equation}\label{e:deltadef}
	\Delta_{f,g} := \{(x,y) \in \RR_+^2 \, : \, f'(x) \neq 0, g'(y) \neq 0,  f'(x) = g'(y)\}.
\end{equation}
Note that since $f(x)$ and $g(y)$ are assumed to be convex decreasing functions, $\Delta_{f,g}$ describes a simple curve in $ [a,b] \times [c,d] \subset \RR_+^2$ if $f'(x) \neq 0$ for all $x \in [a,b]$ and $g'(y) \neq 0$ for all $y \in [c,d]$.

For sets $A, B \subseteq \RR_+^2$ we say that $A$ \emph{ lies Northwest of} $B$ and we write $A \sq B$ if for any $\ul a \in A$ and any $\ul b \in B$ that satisfy $a_1 + a_2 = b_1 + b_2$ we have $a_2 \ge b_2$.

\begin{lemma}\label{lem:var1} If $\gamma_1$ and $\gamma_2$ are paths from $\ul a$ to $\ul b$, and we have either $\gamma_1 \sq \gamma_2 \sq \Delta_{f,g}$ or $\Delta_{f,g} \sq \gamma_2 \sq \gamma_1$, then $w_{f,g} (\gamma_1) \ge w_{f,g}(\gamma_2)$.
\end{lemma}
\proof To start, assume that $\gamma_1 \sq \gamma_2 \sq \Delta_{f,g}$. Let $H$ be the region between $\gamma_1$ and $\gamma_2$:
\[
	H := \{\ul u \, : \, \ul a \le \ul u \le \ul b \text{ and } \gamma_1 \sq \{ \ul u \} \sq \gamma_2\}.
\]
By Green's Theorem in the plane we have
\[
	w_{f,g}(\gamma_1) - w_{f,g}(\gamma_2) = \iint\limits_{H} \left(g'(y) - f'(x) \right) \d x \d y.
\]
Now, since $\gamma_1 \sq \gamma_2 \sq \Delta_{f,g}$ we have $H \sq \Delta_{f,g}$, and since moreover $f$ and $g$ are convex decreasing functions, we have $g'(y) - f'(x) \ge 0$ for all $(x,y) \in H$. It follows that $w_{f,g}(\gamma_1) - w_{f,g}(\gamma_2) \ge 0$.

By the same reasoning we have $w_{f,g}(\gamma_1) - w_{f,g}(\gamma_2) \ge 0$ when $\Delta_{f,g} \sq \gamma_2 \sq~\gamma_1$. \qed
\medskip

\begin{lemma}\label{lem:green} For $\ul a, \ul b \in \Delta_{f,g}$ with $\ul a \le \ul b$, let $\gamma_0 :=\Delta_{f,g} \cap ([a_1, b_1] \times [a_2, b_2])$, then $W_{f,g}(\ul a, \ul b) = w_{f,g}(\gamma_0)$.
\end{lemma}
\proof Suppose by contradiction that $\gamma_1 \neq \gamma_0$ is a minimiser of $W_{f,g}(\ul a, \ul b)$, and $\gamma_0$ is not. Then $\gamma_1$ must intersect $\gamma_0$ in at least two points (counting $\ul a$ and $\ul b$ as intersection points as well). So we can find a set of disjoint curves $\{\eta_i^1\}$ with $\eta_i^1 \subset \gamma_1$ and a set  of disjoint curves $\{\eta_i^0\}$ with $\eta_i^0 \subset \gamma_0$ so that $\gamma_1 \setminus \cup_i \eta_i^1 = \gamma_0 \setminus \cup_i \eta_i^0$ and $\eta_i^1 \sq \eta_i^0 \sq \Delta_{f,g}$ or $\Delta_{f,g} \sq \eta_i^0 \sq \eta_i^1$ for each $i$, and so that $\eta_i^0$ and $\eta_i^1$ have the same end-points. By Lemma~\ref{lem:var1}, replacing the curve $\eta_i^1$ by $\eta_i^0$ in $\gamma_1$ does not increase the value of the line integral. Repeating this procedure for each such interval, we end up replacing the minimiser $\gamma_1$ by $\gamma_0$ without increasing the value of the integral, contradicting the assumption that $\gamma_0$ was not a minimiser. \qed
\medskip
 
Given a set of points $\{\ul a_{\sss (i)}\}$, with $\ul a_{\sss (i)} \in \RR_+^2$, we write $\ul a_{\sss (1)} \to \ul a_{\sss (2)} \to  \dotsm \to \ul a_{\sss (n)}$ for the path that linearly interpolates between successive points $\ul a_{\sss (i)}$ and $\ul a_{\sss (i+1)}$. Given a path $\gamma$ and two points $\ul a, \ul b \in \gamma$, we write $\ul a \stackrel{\gamma}{\to} \ul b$ for the part of $\gamma$ between $\ul a$ and $\ul b$.

  \begin{lemma}\label{lem:easypath} If $g(y)=c$ for some constant $c \ge 0$ and $f(x)$ is a positive, monotone decreasing function, then for $\ul a \le \ul b$,
 \[
 	W_{f,g}(\ul a, \ul b) = c (b_1-a_1) + f(b_1) (b_2 -a_2),
\]
and the path from $\ul a$ to $\ul b$ that minimises $W_{f,g}(\ul a, \ul b)$ is $\ul a \to (b_1, a_2) \to \ul b$.
\end{lemma}
\proof This follows directly from the definition of $W_{f,g}$ and the assumptions on $f$ and~$g$.\qed
\medskip

Recall the definitions of $\psi(y)$ and $\phi(x)$ from \eqref{e:gdef} and \eqref{e:fdef}, and the definition of $\Delta_{f,g}$ in \eqref{e:deltadef}.
Observe that
\[
	\psi'(x) = \frac{1}{x} \indi_{\{\frac{3 \xi^2}{p} \le x \le \frac{1}{p^2}\}} \qquad \text{ and } \qquad \phi'(y) = -3p \e^{-3py} \indi_{\{y > \frac{4}{p} \log \log \frac1p\}},
\]
so $\psi'(x) = \phi'(y)$ is solved by
\[
	x(y) = \frac{e^{3py}}{3p} \qquad \text{ and } \qquad y(x) = \frac{1}{3p} \log(3px)
\]
when both $\psi(x) \neq 0$ and $\phi(y) \neq 0$. Observe that
\[\begin{split}
	y\left(\frac{3 \xi^2}{p}\right) &= \frac{1}{3p} \log(9 \xi^2) < \frac{4}{p} \log\log \frac1p,\\
	y\left(\frac{1}{p^2}\right) &= \frac{1}{3p} \log \frac3p,\\
	x\left(\frac{4}{p} \log \log \frac 1p\right) &= \frac{1}{3p} \log^{12} \frac1p.
\end{split}
\]
We can thus write
\[\begin{split}
	\Delta_{\psi,\phi}  & = \left\{(x,y) \in \RR_+^2 \, : \, \psi'(x) \neq 0,  \phi'(y) \neq 0,  \psi'(x) = \phi'(y) \right\} \\ 
		& =  \left\{\left(\frac{\e^{3py}}{3 p}, y\right) \, : \, y \in \left(\frac{4}{p} \log\log \frac1p, \frac{1}{3p} \log \frac3p \right] \right\}.
\end{split}\]
 The leftmost and rightmost points of $\Delta_{\psi,\phi}$ are given by
 \begin{equation}\label{e:uvdef}
 	\ul u = \left(\frac{1}{3  p} \log^{12} \frac1p, \frac{4}{p}  \log\log\frac1p \right) \quad \text{ and } \quad \ul v = \left(\frac{1}{p^2}, \frac{1}{3p} \log \frac3p\right).
\end{equation}
\begin{lemma}\label{lem:overshoot}  Let $\ul u $ and $\ul v$ be as in \eqref{e:uvdef}, and let $\ul a$ be such that $\ul a \le \ul u$, and let $\ul b$ be such that $\ul v \le \ul b$. Then
\[
	W_{\psi,\phi}(\ul a, \ul b) = W_{\psi,\phi}(\ul a, \ul u) + W_{\psi,\phi}(\ul u, \ul v) + W_{\psi, \phi}(\ul v, \ul b).
\]
\end{lemma}
\proof 
By Lemma~\ref{lem:12}, the right-hand side is an upper bound on $W_{\psi,\phi}(\ul a, \ul b)$. It remains to prove that it is also a lower bound.

Since $\psi$ and $\phi$ are decreasing, positive, continuous functions, any path that minimises $W_{\psi,\phi}(\ul a, \ul b)$ must be a coordinate-wise increasing path.
Fix a coordinate-wise increasing path $\gamma \subset \RR_+^2$ from $\ul a$ to $\ul b$. 
Then either
\begin{enumerate}
	\item $\gamma \cap \Delta_{\psi,\phi} \neq \emptyset,$ or
	\item $\gamma \cap \Delta_{\psi,\phi} = \emptyset.$
\end{enumerate}

\begin{figure}
	\includegraphics[width = .49\textwidth]{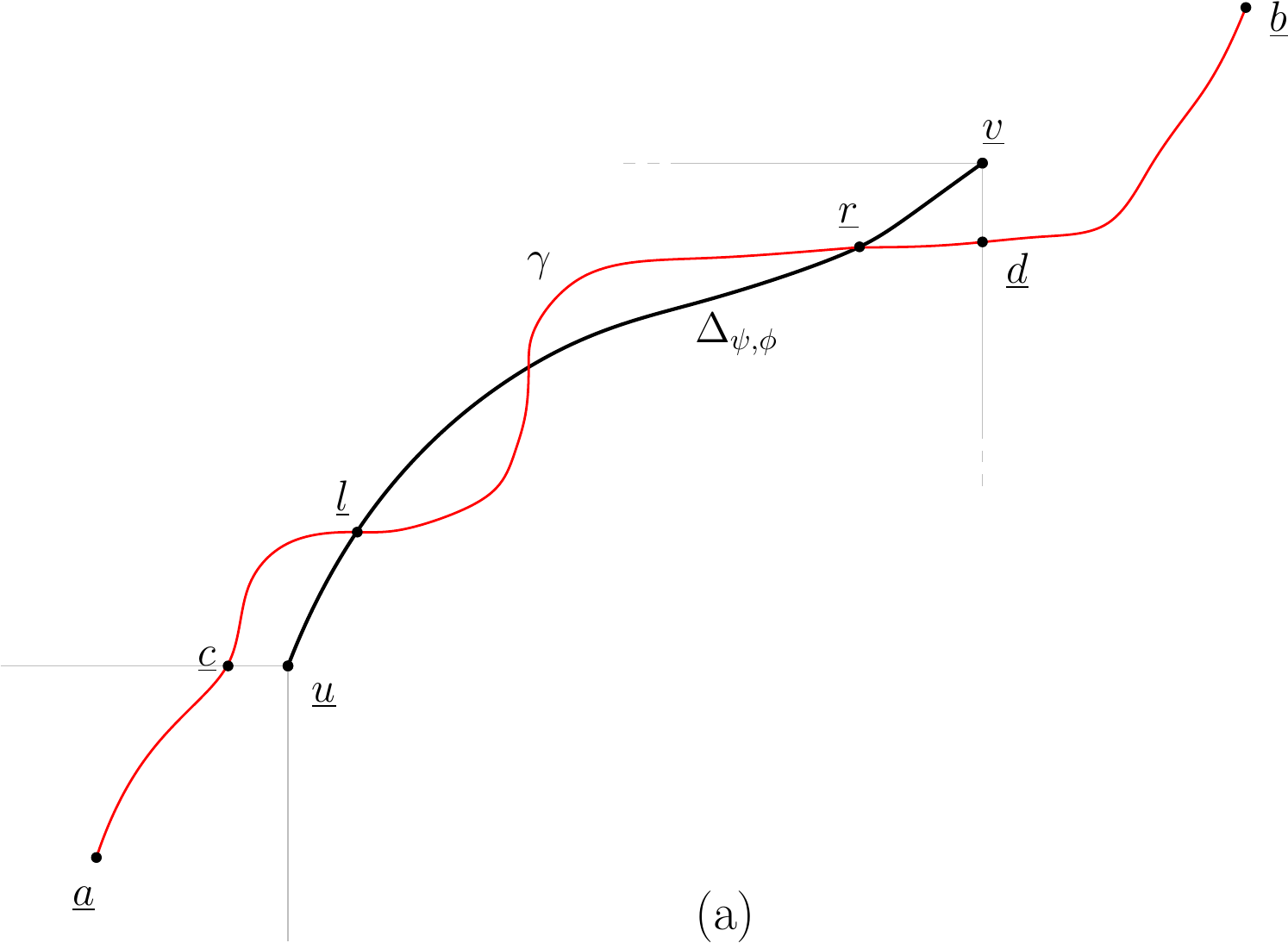}
	\includegraphics[width = .49\textwidth]{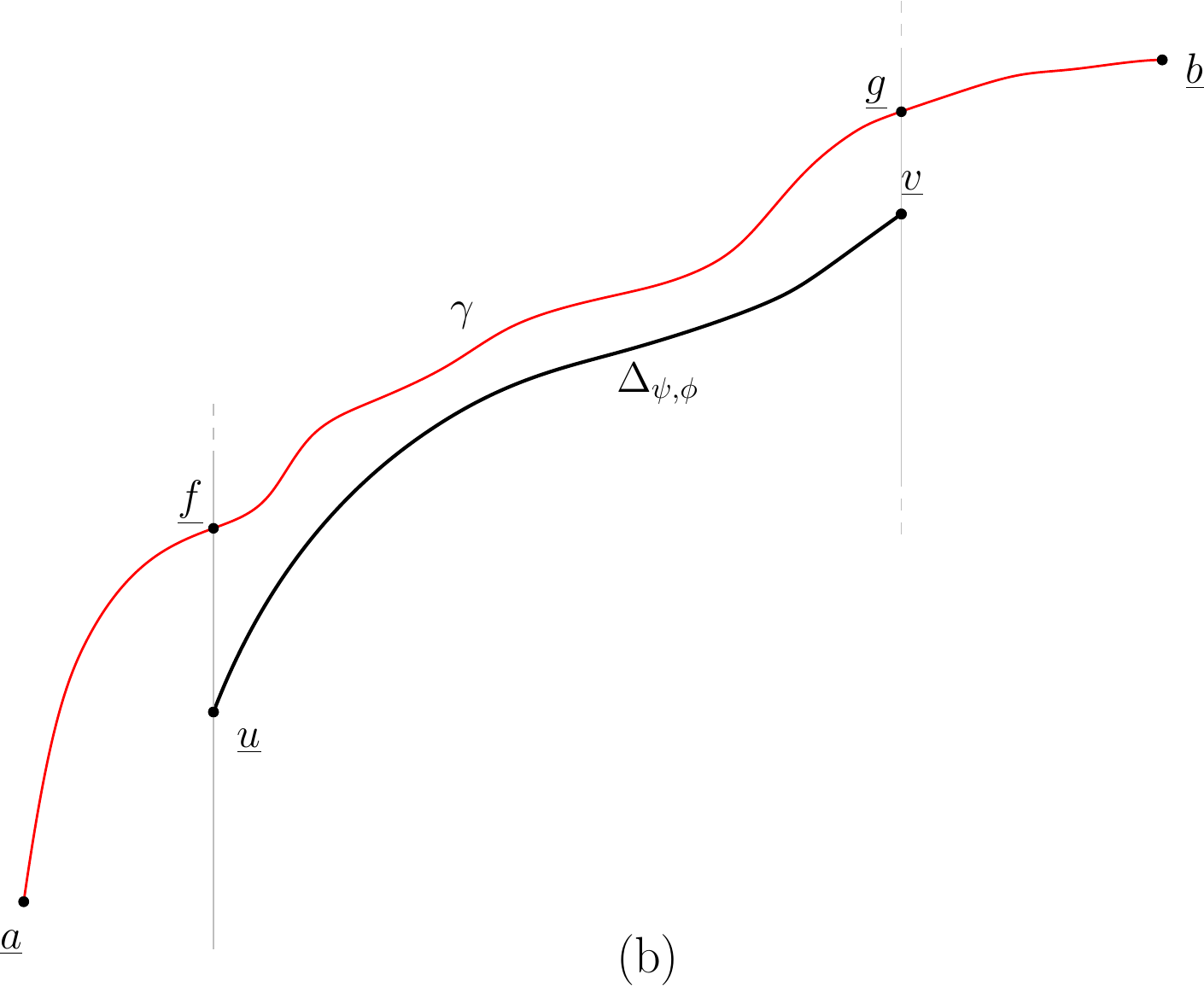}
	\caption{\label{fig:Lem56} The partitioning of $\gamma$ used in Lemma \ref{lem:overshoot} for the case $\gamma \cap \Delta_{\psi,\phi} \neq \emptyset$ in (a), and for the case $\gamma \cap \Delta_{\psi,\phi} = \emptyset$ in (b).}
\end{figure}

Consider first case (a). Write $\ul c \in \gamma$ for the first point in $\gamma$ such that either $c_1 = u_1$ or $c_2 = u_2$ and $\ul d \in \gamma$ for the first point in $\gamma$ such that either $d_1 = v_1$ or $d_2 = v_2$. Write $\ul l$ and $\ul r$ for the first and last point along $\gamma$ where $\gamma$ and $\Delta_{\psi, \phi}$ intersect. Since $\gamma$ is coordinate-wise increasing we have $\ul a \le \ul c \le \ul l \le \ul r \le \ul d \le \ul b$. See Figure~\ref{fig:Lem56}(a). We split the integral along $\gamma$ into five parts:
\[
	w_{\psi,\phi}(\gamma) = w_{\psi,\phi}(\ul a \stackrel{\gamma}{\to} \ul c) + w_{\psi,\phi}(\ul c \stackrel{\gamma}{\to} \ul l) + w_{\psi,\phi}(\ul l \stackrel{\gamma}{\to} \ul r) + w_{\psi,\phi}(\ul r \stackrel{\gamma}{\to} \ul d) + w_{\psi,\phi}(\ul d \stackrel{\gamma}{\to} \ul b).
\]
Using Lemma \ref{lem:green} we split the minimising integral from $\ul u$ to $\ul v$ into three parts:
\[
	W_{\psi, \phi}(\ul u, \ul v) = w_{\psi,\phi} \big(\ul u \stackrel{\Delta_{\psi,\phi}}{\longrightarrow} \ul l \big) + w_{\psi,\phi}\big(\ul l \stackrel{\Delta_{\psi,\phi}}{\longrightarrow} \ul r\big) + w_{\psi,\phi}\big(\ul r \stackrel{\Delta_{\psi,\phi}}{\longrightarrow} \ul v\big).
\]
By Lemma \ref{lem:easypath},
\[
	w_{\psi, \phi}(\ul a \stackrel{\gamma}{\to} \ul c) \ge w_{\psi, \phi}(\ul a \to (c_1, a_2) \to \ul c).
\]
By Lemma \ref{lem:var1} and the fact that either 
\[
	\ul c \stackrel{\gamma}{\to} \ul l \sq \ul c \to \ul u \stackrel{\Delta_{\psi,\phi}}{\longrightarrow} \ul l  \sq \Delta_{\psi,\phi} \qquad \text{ or } \qquad \Delta_{\psi,\phi} \sq \ul c \to \ul u \stackrel{\Delta_{\psi,\phi}}{\longrightarrow} \ul l \sq \ul c \stackrel{\gamma}{\to} \ul l,
\]
we have
\[
	w_{\psi, \phi}(\ul c \stackrel{\gamma}{\to} \ul l) \ge w_{\psi, \phi}(\ul c \to \ul u) + w_{\psi, \phi}\big(\ul u \stackrel{\Delta_{\psi,\phi}}{\longrightarrow} \ul l\big).
\]
By Lemma \ref{lem:12},
\[
	w_{\psi, \phi}(\ul a \to (c_1, a_2) \to \ul c) + w_{\psi, \phi}(\ul c \to \ul u) = w_{\psi, \phi}(\ul a \to (c_1, a_2) \to \ul c \to \ul u) \ge W_{\psi,\phi}(\ul a, \ul u).
\]
Moreover, by Lemma \ref{lem:green},
\[
	w_{\psi, \phi}(\ul l \stackrel{\gamma}{\to} \ul r) \ge w_{\psi, \phi}\big(\ul l \stackrel{\Delta_{\psi,\phi}}{\longrightarrow} \ul r\big).
\]
By Lemma \ref{lem:var1} and the fact that either 
\[
	\ul r \stackrel{\gamma}{\to} \ul d \sq \ul r \stackrel{\Delta_{\psi,\phi}}{\longrightarrow} \ul v \to \ul d \sq \Delta_{\psi,\phi} \qquad \text{ or } \qquad \Delta_{\psi,\phi} \sq \ul r \stackrel{\Delta_{\psi,\phi}}{\longrightarrow}\ul v \to \ul d \sq \ul r \stackrel{\gamma}{\to} \ul d,
\]
we have
\[
	w_{\psi, \phi}(\ul r \stackrel{\gamma}{\to} \ul d) \ge w_{\psi, \phi}\big(\ul r \stackrel{\Delta_{\psi,\phi}}{\longrightarrow} \ul v \big) + w_{\psi, \phi}(\ul v \to \ul d).
\]
And finally, since $\ul v \to \ul d \stackrel{\gamma}{\to} \ul b$ is a path from $\ul v$ to $\ul b$,
\[
	W_{\psi, \phi}(\ul v, \ul b) \le w_{\psi, \phi}(\ul v \to \ul d) + w_{\psi, \phi}(\ul d \stackrel{\gamma}{\to} \ul b).
\]

Combining the above inequalities we obtain
\begin{equation}\label{e:gammaintersect}\begin{split}
	w_{\psi,\phi}(\gamma) &\ge W_{\psi, \phi}(\ul a , \ul u) + w_{\psi, \phi}\big(\ul u \stackrel{\Delta_{\psi,\phi}}{\longrightarrow} \ul l\big) + w_{\psi, \phi}\big(\ul l \stackrel{\Delta_{\psi,\phi}}{\longrightarrow} \ul r\big)\\
	& \qquad +  w_{\psi, \phi}\big(\ul r \stackrel{\Delta_{\psi,\phi}}{\longrightarrow} \ul v\big) + w_{\psi, \phi}(\ul v \to \ul d) + W_{\psi, \phi}(\ul v, \ul b) - w_{\psi, \phi}(\ul v \to \ul d)\\
	& \ge W_{\psi,\phi}(\ul a, \ul u) + W_{\psi,\phi}(\ul u, \ul v)  + W_{\psi,\phi}(\ul v, \ul b).
\end{split}
\end{equation}
\medskip

Now we consider case (b), that $\gamma \cap \Delta_{\psi,\phi} = \emptyset$. Let $\ul f$ be the first point on $\gamma$ such that $f_1 = u_1$, and let $\ul g$ be the first point on $\gamma$ such that $g_1 = v_1$. See Figure~\ref{fig:Lem56}(b). We divide the integral along $\gamma$ into three parts:
\[
	w_{\psi,\phi}(\gamma) = w_{\psi, \phi}(\ul a \stackrel{\gamma}{\to} \ul f)+ w_{\psi, \phi}(\ul f \stackrel{\gamma}{\to} \ul g) + w_{\psi, \phi}(\ul g \stackrel{\gamma}{\to} \ul b).
\]
By Lemma \ref{lem:easypath},
\[\begin{split}
	w_{\psi,\phi}(\ul a \stackrel{\gamma}{\to} \ul f) &\ge w_{\psi,\phi}(\ul a \to (f_1,a_2) \to \ul f) = w_{\psi,\phi}(\ul a \to (u_1, a_2 ) \to \ul u) + w_{\psi,\phi}(\ul u \to \ul f)\\
	&= W_{\psi,\phi}(\ul a, \ul u) + w_{\psi,\phi}(\ul u \to \ul f).
\end{split}\]
Since $\gamma \cap \Delta_{\psi,\phi} = \emptyset$, either $f_2 \ge u_2$ and $g_2 \ge v_2$ or $f_2 \le u_2$ and $g_2 \le v_2$, so that either
\[
	\ul u \to \ul f \stackrel{\gamma}{\to} \ul g \to \ul v \sq \Delta_{\psi,\phi} \qquad \text{ or } \qquad \Delta_{\psi,\phi}
 \sq \ul u \to \ul f \stackrel{\gamma}{\to} \ul g \to \ul v .
\]
It thus follows by Lemmas \ref{lem:var1} and \ref{lem:green} that 
\[
	w_{\psi,\phi}(\ul f \stackrel{\gamma}{\to} \ul g) \ge W_{\psi,\phi}(\ul u, \ul v) - w_{\psi,\phi}(\ul u \to \ul f) - w_{\psi, \phi}(\ul g \to \ul v).
\]
 Finally, since $\ul v \to \ul g \stackrel{\gamma}{\to} \ul b$ is a path from $\ul v$ to $\ul b$, 
 \[
 	W_{\psi,\phi}(\ul v, \ul b) \le w_{\psi,\phi}(\ul v \to \ul g) + w_{\psi,\phi}(\ul g \stackrel{\gamma}{\to} \ul b).
\]

Combining the above inequalities we obtain
\begin{equation}\label{e:gammanointersect}\begin{split}
	w_{\psi,\phi}(\gamma) &\ge W_{\psi,\phi}(\ul a, \ul u) + w_{\psi,\phi}(\ul u \to \ul f) + W_{\psi,\phi}(\ul u, \ul v) - w_{\psi,\phi}(\ul u \to \ul f)\\
	& \qquad  - w_{\psi, \phi}(\ul g \to \ul v) + W_{\psi,\phi}(\ul v, \ul b) - w_{\psi,\phi}(\ul v \to \ul g)\\
	& = W_{\psi,\phi}(\ul a, \ul u) + W_{\psi,\phi}(\ul u, \ul v)  + W_{\psi,\phi}(\ul v, \ul b).
\end{split}
\end{equation}

From \eqref{e:gammaintersect} and \eqref{e:gammanointersect} we conclude that the lower bound $w_{\psi,\phi}(\gamma) \ge W_{\psi,\phi}(\ul a, \ul u) + W_{\psi,\phi}(\ul u, \ul v) + W_{\psi,\phi}(\ul v, \ul b)$ holds uniformly for any path $\gamma$, and therefore, it also holds for the infimum, completing the proof.
\qed
\medskip

The following lemma now states the crucial bound:
 \begin{lemma}\label{lem:Wbd2} 
 Let $\ul u$ and $\ul v$ be as in \eqref{e:uvdef}, and let  $\ul a \le \ul u$, $a_2 = o(1/p)$, and $\ul b \ge \ul v$. Then
 \[
 	W_{\psi,\phi}\left(\ul a, \ul b \right) \ge \frac{1}{6p} \log^2 \frac1p - (1+o(1))\frac{1}{3p} \log\left(\frac{8}{3 \e}\right) \log \frac1p.
\]
\end{lemma}
 \proof
Since we assumed $\ul a \le \ul u$ and $\ul b \ge \ul v$, Lemma \ref{lem:overshoot} gives
\begin{equation}\label{e:splitW}
	W_{\psi,\phi}(\ul a, \ul b) \ge W_{\psi,\phi}(\ul a, \ul u) + W_{\psi,\phi}(\ul u, \ul v).
\end{equation}

Recall that we have set $\xi = \lceil\log^2 \frac1p \rceil$ and $\delta_\xi = 1-2/\xi$.
We use Lemma \ref{lem:easypath}, that $a_1 \le u_1$, and that $a_2 = o(1/p)$ to bound
\begin{equation}\label{e:Wstart}
\begin{split}
	W_{\psi,\phi}(\ul a, \ul u) \ge & - \left( \log(8p^2 u_1 + 8p) +1/\xi\right)  \left(u_2- a_2\right)\\
	=& \frac{4}{p} \log \frac1p \log\log\frac1p - o\left(\frac1p \log \frac1p\right).
\end{split}
\end{equation}

Now we bound $W_{\psi,\phi}(\ul u, \ul v)$.
It follows by Lemma \ref{lem:green} that
\begin{equation}\label{e:integral}\begin{split}
	W_{\psi,\phi}\left(\ul u, \ul v\right)  & = \int\limits_{\frac{4}{p}\log\log \frac1p}^{\frac{1}{3p} \log \frac3p} \psi(x(y)) \d y + \phi( y) x'(y) \d y \\
		&\ge  \int\limits_{\frac{4}{p}\log\log \frac1p}^{\frac{1}{3p} \log \frac{1}{p}} \left(- \log\left(\frac{8p}{3} \e^{3py} \right)  - 1/\xi +   \e^{-3py} \cdot \e^{3py}\right) \d y\\
		&\ge  \int\limits_{\frac{4}{p}\log\log \frac1p}^{\frac{1}{3p} \log \frac{1}{p}}\left(1 - 3/\xi - \log\left(\frac{8p}{3}\right) - 3py \right)  \d y.
\end{split}
\end{equation}
The integral on the right-hand side evaluates to
\begin{multline*}
	 \left(1 - 3/\xi -  \log \frac{8p}{3} \right)\Bigr [ y \,\Bigl ]_{\frac{4}{p}\log\log \frac1p}^{\frac{1}{3p} \log \frac{1}{p}} - \left[ \frac{3 p y^2}{2} \right]_{\frac{4}{p}\log\log \frac1p}^{\frac{1}{3p} \log \frac{1}{p}}\\
	=  \frac{1}{6p} \log^2 \frac1p - \frac{4}{p} \log\frac1p \log\log \frac1p - (1+o(1))\frac{1}{3p} \log \frac{8}{3 \e} \log \frac1p.
\end{multline*}
(Observe that the first term in the first integral in \eqref{e:integral} thus gives a complementary bound to \eqref{e:prod2}, while the second term is complementary to \eqref{e:prod1}.)
It follows that
\begin{equation}\label{e:wsum}
	W_{\psi,\phi}\left(\ul u, \ul v \right)  \ge \frac{1}{6p} \log^2 \frac1p -  \frac{4}{p} \log\frac1p \log\log \frac1p - (1+o(1))\frac{1}{3p} \log \frac{8}{3\e} \log \frac1p .
\end{equation}

Substituting \eqref{e:Wstart} and \eqref{e:wsum} into \eqref{e:splitW} completes the proof. \qed
\medskip
 
 \proof[Proof of Lemma \ref{lem:Wbd}] Given a sequence of increasing rectangles $(R_n)_{n=0}^N$, let $(x_n, y_n) \in \RR_+^2$ denote the dimensions of $R_n$. Construct the path $\gamma \subset \RR_+^2$ by linearly interpolating between successive points $(x_n, y_n)$, i.e.,
 \[
 	\gamma = (x_0, y_0) \to (x_1, y_1) \to \dotsm \to (x_N, y_N).
\]
Recall the definition of $U^p(R,R')$ from \eqref{e:Wdef} and recall that $\delta_\xi  = 1 - \log^{-2} \frac1p$. It follows from Lemma \ref{lem:12} that 
\[
	\sum_{n=0}^{N-1} U^p (R_n, R_{n+1}) 
	\ge \delta_\xi W_{\psi, \phi}\left((x_0, y_0), (x_{N}, y_{N})\right),   
\]
and it follows from Lemma \ref{lem:Wbd2} that the right-hand side is bounded from below by
\[
	 \frac{1}{6p} \log^2 \frac1p - (1+o(1))\frac{1}{3p} \log\left(\frac{8}{3 \e}\right) \log \frac1p,
\]
completing the proof. \qed

\section{The critical probability: proof of Theorem \ref{thm:pc}}\label{sec:pc}
We start with the upper bound. From \cite{DumEnt13} we know that if $L \ll \e^{p^{-1+\vep}}$ for any $\vep >0$, then $\Pp([L]^2$ is IF$) = o(1)$, so we assume that $L \ge \e^{p^{-1+\vep}}$. 
Let $m=p^{-5}$. The probability that $[L]^2$ is internally filled is bounded from below by the probability that $[L]^2$ contains exactly one internally filled translate of $[m]^2$, and that $\{[m]^2  \growto [L]^2\}$ occurs.
Indeed, let $R_{\sss (i,j)} = [(i-1)m+1, im] \times [(j-1)m+1, jm]$, and let
\[
	A_{(i,j)} := \{R_{\sss (i,j)} \text{ is IF}\} \cap \{\forall (i',j') \neq (i,j): R_{\sss (i',j')} \text{ is not IF}\} \cap \{R_{\sss (i,j)} \growto [L]^2\}.
\]
Then
\[
	\{[L]^2 \text{ is IF}\} \supseteq \bigsqcup_{i,j = 1}^{\lfloor L / m\rfloor} A_{\sss (i,j)},
\]
so that
\begin{multline}\label{e:blablabla}
	\Pp([L]^2 \text{ is IF}) 
	\ge \sum_{i,j=1}^{\lfloor L/m\rfloor} \bigg( \Pp (R_{\sss (i,j)}\text{ is IF}) - \big(1-\Pp(R_{\sss (i,j)} \growto [L]^2)\big) \\
	- \sum_{\substack{i',j'=1:\\
	i' \neq i , j' \neq j}}^{\lfloor L / m\rfloor} \Pp(R_{\sss (i,j)} \text{ and } R_{\sss (i',j')} \text{ are IF})\bigg).
	\end{multline}

To bound $\Pp(R_{\sss (i,j)} \growto [L]^2)$, observe that if every horizontal and vertical line segment of length $p^{-5}$ intersecting $[L]^2$ contains a pair of adjacent infected sites, then it must be the case that $\{R_{\sss (i,j)} \growto [L]^2\}$ occurs. We can bound the probability of this event from below by
\[
	(1-(1-p^2)^{\frac{m}{2}})^{4 L^2} = \exp\big(4 L^2 \e^{-\frac{1}{2 p^{3}}(1+o(1))}\big) \ge \exp( c \e^{-c' p^{-3}}),
\]
for some $c, c'>0$ and $p$ sufficiently small. Note that by Proposition~\ref{upper}, the right-hand side is $1-o(\Pp([m]^2 \text{ is IF}))$.
Inserting this bound into \eqref{e:blablabla} and summing over the indices, we obtain
\[\begin{split}
	\Pp([L]^2 \text{ is IF}) &\ge \frac{L^2}{m^2} \bigg(\Pp([m]^2 \text{ is IF}) - \big(1-\exp( c \e^{-c' p^{-3}})\big) - \frac{L^2}{m^2} \Pp([m]^2 \text{ is IF})^2 \bigg)\\
	& \ge \frac{L^2}{2 m^2} \Pp([m]^2 \text{ is IF})\Big(1 - \frac{2 L^2}{m^2} \Pp([m]^2 \text{ is IF})\Big),
\end{split}\]
where the second inequality holds for $p$ sufficiently small.
Taking minus the logarithm on both sides and applying the above bound and Proposition~\ref{upper} again we obtain the inequality
\begin{multline*}
	 -\log(\Pp([L]^2 \text{ is IF})) \le 2C_1 \frac1p \log^2 \frac1p - 2C_2 \frac1p \log \frac1p  - 2 \log L + 10 \log \frac1p + \log 2 \\
	 - \log\bigg(1- \frac{2L^2}{p^{10}} \exp\Big(-2C_1 \frac1p \log^2 \frac1p + (2C_2-o(1)) \frac1p \log \frac1p\Big)\bigg).
\end{multline*}
Observe that if the right-hand side tends to $0$ from above when we let $p \to 0$, then $\Pp \left([L]^2 \text{ is IF}\right) \to 1$. 
To make the right-hand side vanish we will fix $L$ to be equal to $\Lambda$, where 
\[
	\Lambda = \Lambda(p) := \exp\left( \frac{C_1}{p}  \log^2 \frac 1p  - \frac {C_2 -\eta_u}{p} \log \frac 1p  \right),
\]
where $\eta_u :=  p (10 + \log 2 \log^{-1} \frac1p ) = o(1)$. Indeed, with this choice all the leading order terms cancel, and we obtain
\[
	-\log(\Pp([L]^2 \text{ is IF})) \le -\log\bigg(1- 4 \exp\Big(-o\Big(\frac1p \log \frac 1p \Big)\Big)\bigg) \to 0 \text{ as } p \to 0,
\]
as desired.

Now we invert $\Lambda(p)$ to find $p_u = p_u (L)$, an asymptotically minimal sequence in $p$ such that 
\[
	\Pr_{p_u(L)} ([L]^2 \text{ is IF}) \to 1 \qquad \text{ as } \qquad L  \to \infty.
\]
We can express $p_u$ and $1/p_u$ in terms of $\Lambda$:
\begin{equation}\label{eq-pLambda}
	p_u = \frac{C_1\log^2 \frac {1}{p_u}  - (C_2-\eta_{u}) \log \frac {1}{p_u}}{ \log \Lambda},
\end{equation}
and
\begin{equation}\label{eq-invpLambda}
	\frac {1}{p_u} =  \frac{ \log \Lambda}{C_1  \log^2 \frac {1}{p_u}  - (C_2- \eta_{u}) \log \frac {1}{p_u}}.
\end{equation}
Substituting \eqref{eq-invpLambda} into \eqref{eq-pLambda}, we get
\begin{multline*}
	p_u = \frac{1}{\log \Lambda} \Biggl(C_1 \left(\left(\log\log \Lambda \right) - \log\left(C_1 \log^2 \frac {1}{p_u} - (C_2-\eta_{u}) \log \frac {1}{p_u} \right) \right)^2 \\
	- (C_2-\eta_{u}) \left(\left(\log\log \Lambda \right) - \log \left(C_1 \log^2 \frac {1}{p_u}  - (C_2-\eta_{u}) \log \frac {1}{p_u} \right) \right)\Biggr).
\end{multline*}
For sufficiently large $\Lambda$ (and hence for small $p_u$) and $\delta>0$, the following inequalities hold:
\[
	 \frac 1{p_u}  \le  \log \Lambda  \le  \frac{1}{p_u^{1 + \delta}}, \qquad \text{ and } \qquad
	 \log \frac 1{p_u}  \le  \log\log \Lambda  \le  (1+ \delta) \log \frac 1{p_u}.
\]
Whence we obtain the asymptotic formula
\[
	p_u = \frac{C_1 (\log\log \Lambda)^2}{\log \Lambda} - \frac{4 C_1 \log\log \Lambda \log\log\log \Lambda}{\log \Lambda} 
	+ \frac{\left(C_2 + 2 C_1 \log C_1 + \delta \right)\log\log \Lambda}{\log \Lambda},
\]
for $\delta >0$, giving the upper bound in Theorem \ref{thm:pc}.
\medskip

Now we prove the lower bound. Again let $m=p^{-5}$ and let $n = \frac1p \log \frac1p$ and let $L \gg m,n$.  Let $\R$ denote the set of all rectangles $R \subset [L]^2$ with dimensions $(x,y)$ such that $m/3 \le x \le m$ and $n/3 \le y \le n$. It is a straightforward consequence of the proof of \cite[Lemma 3.7]{DumEnt13} that if $[L]^2$ is internally filled, then there must exist a rectangle $R \in \R$ such that $\{R \growto [L]^2\}$ occurs.\footnote{This follows if we stop the algorithm in the proof of \cite[Lemma 3.7]{DumEnt13} when a set $\S$ such that $\langle \S \rangle \subset \R$ is first constructed, and by observing that if $[L]^2$ is internally filled, then we must construct such a set $\S$.} The number of rectangles in $\R$ is bounded by $ mn L^2$, so by Proposition~\ref{lower},
\[ \begin{split}
	\Pp \left([L]^2 \text{ is IF}\right)& = \Pp \Big(\bigcup_{R \in \R} \left\{ R \text{ is IF} \right\} \cap \left\{ R \growto [L]^2 \right\} \Big)\\
			&\le  m n L^2  \exp\Big(-2C_1\frac{1}{p} \log^2 \frac1p + (2C_2 + \zeta)\frac{1}{p} \log \frac1p \Big) 
\end{split}
\]
for any $\zeta > 0$.
Taking the logarithm of both sides gives
\[
		 \log \Pp \left([L]^2 \text{ is IF}\right) 
					 \le 2 \log L  - 2C_1 \frac1p \log^2 \frac 1p + 2(C_2 + \eta_\ell) \frac1p\log \frac 1p,
\]
where $\eta_\ell := p \log(mn) + \zeta$.
Observe that if the upper bound tends to $- \infty$ as $p \to 0$, then $\Pp \left([L]^2 \text{ is IF}\right) \to 0$.
To minimise the right-hand side, we will fix $L$ to be equal to $\lambda$, where
\[
	\lambda  = \lambda(p) := \exp\left( \frac{C_1}{p}  \log^2 \frac 1p  - \frac {C_2 +\eta_\ell}{p} \log \frac 1p  \right).
\]
Again, we invert $\lambda(p)$, now to find $p_\ell = p_\ell(L)$, an asymptotically maximal sequence in $p$ such that 
\[
	\Pr_{p_\ell(L)} ([L]^2 \text{ is IF}) \to 0 \qquad \text{ as } \qquad L \to \infty.
\] 
Using the same steps as we used in the proof of the upper bound, we can now determine that $p_\ell$ satisfies
\[
	p_\ell = \frac{C_1 (\log\log \lambda)^2}{\log \lambda} - \frac{4 C_1\log\log \lambda \log\log\log \lambda}{ \log \lambda} 
	+ \frac{\left(C_2 + 2 C_1 \log C_1 - \delta \right)\log\log \lambda}{\log \lambda}
\]
for some $\delta >0$ that can be chosen arbitrarily small but depends on $\zeta$.
This concludes the proof of Theorem \ref{thm:pc}.
\qed

\begin{small}
\bibliographystyle{abbrv}
\bibliography{bootstrap}
\end{small}
\end{document}